\documentclass[nohyperref]{article}

\usepackage{microtype}
\usepackage{graphicx}
\usepackage{subfigure}
\usepackage{booktabs} %

\usepackage{hyperref}

\usepackage[accepted]{icml2022}

\usepackage{amsmath}
\usepackage{amssymb}
\usepackage{mathtools}
\usepackage{amsthm}

\usepackage[capitalize,noabbrev]{cleveref}

\usepackage[textsize=tiny]{todonotes}

\def\d{\textnormal{d}}

\def\dist{\textnormal{dist}}

\def\conv{\mathop{\textnormal{Co}}}

\def\dom{\mathop{\textnormal{dom}}}

\def\vvec{\mathop{\textnormal{Vec}}}

\usepackage[utf8]{inputenc} %
\usepackage[T1]{fontenc}    %

\usepackage{tcolorbox}

\usepackage{microtype}
\usepackage{graphicx}
\usepackage{booktabs} %

\usepackage{amsmath,amssymb,amsthm}
\usepackage{bbold}
\usepackage[sans]{dsfont}

\usepackage{cleveref}

\usepackage{pifont}

\usepackage{algorithm,algorithmic}

\usepackage{thmtools}
\usepackage{thm-restate}

\usepackage{multirow}

\usepackage{tabularx}

\usepackage{nicematrix}

\def\O#1{\text{\ding{\the\numexpr#1+171}}}

\declaretheorem[name=Theorem,within=section]{theorem}

\declaretheorem[name=Fact,sibling=theorem]{fact}

\declaretheorem[name=Remark,sibling=theorem]{remark}
\declaretheorem[name=Definition,sibling=theorem]{definition}
\declaretheorem[name=Claim,sibling=theorem]{claim}

\declaretheorem[name=Assumption,sibling=theorem]{assumption}

\usepackage{amsmath,amssymb}
\usepackage{tikz}
\usetikzlibrary{shadows}
\usetikzlibrary{arrows,positioning} 
\usetikzlibrary{shapes,chains}
\usetikzlibrary{patterns}

\usepackage{contour}
\usepackage[normalem]{ulem}

\contourlength{0.8pt}

\numberwithin{equation}{section}

\usepackage{bm}

\usepackage{wrapfig}

\definecolor{amethyst}{rgb}{0.6, 0.4, 0.8}

\usepackage{mathbbol}

\newcommand{\sg}{g} %
\newcommand{\so}[1]{\mathbb{O}_{\textup{s}}(#1)} %

\newcommand{\pf}[1]{\nabla f (#1)}
\newcommand{\norm}[1]{\left\lVert#1\right\rVert}
\newcommand{\innr}[1]{\left\langle#1\right\rangle}
\newcommand{\E}[1]{\mathbb{E} \left[#1\right]}
\newcommand{\R}{\mathbb{R}}
\newcommand{\mleq}[1]{\overset{\mathclap{(#1)}}{\leq}}

\newcommand{\mar}[1]{\left(#1\right)}

\icmltitlerunning{Practical Optimization of Lipschitz Functions with Finite-Time Complexity}

\begin{document}

\twocolumn[
\icmltitle{On the Finite-Time Complexity and Practical Computation of \\ Approximate Stationarity Concepts of Lipschitz Functions}

\begin{icmlauthorlist}
\icmlauthor{Lai Tian}{seem}
\icmlauthor{Kaiwen Zhou}{cse}
\icmlauthor{Anthony Man-Cho So}{seem}
\end{icmlauthorlist}

\icmlaffiliation{seem}{Department of Systems Engineering \& Engineering Management, The Chinese University of Hong Kong, Shatin, NT, Hong Kong}
\icmlaffiliation{cse}{Department of Computer Science and Engineering, The Chinese University of Hong Kong, Shatin, NT, Hong Kong}

\icmlcorrespondingauthor{Anthony Man-Cho So}{manchoso@se.cuhk.edu.hk}

\icmlkeywords{Machine Learning, ICML}

\vskip 0.3in
]

\printAffiliationsAndNotice{}  %

\begin{abstract}
We report a practical finite-time algorithmic scheme to compute approximately stationary points for nonconvex nonsmooth Lipschitz functions.
In particular, we are interested in two kinds of approximate stationarity notions for nonconvex nonsmooth problems, i.e., Goldstein approximate stationarity (GAS) and near-approximate stationarity (NAS). For GAS, our scheme removes the unrealistic subgradient selection oracle assumption in \citep[Assumption 1]{zhang2020complexity} and computes GAS with the same finite-time complexity. For NAS, \citet{davis2019stochastic} showed that $\rho$-weakly convex functions admit finite-time computation, while \citet{tian2021hardness} provided the matching impossibility results of dimension-free finite-time complexity for first-order methods. Complement to these developments, in this paper, we isolate a new class of functions that could be Clarke irregular (and thus not weakly convex anymore) and show that our new algorithmic scheme can compute NAS points for functions in that class within finite time. 
To demonstrate the wide applicability of our new theoretical framework, we show that $\rho$-margin SVM, $1$-layer, and $2$-layer ReLU neural networks, all being Clarke irregular, satisfy our new conditions.
\end{abstract}

\section{Introduction}

In this paper, we consider the following general optimization problem for an $L$-Lipschitz function $f:\mathbb{R}^d\rightarrow \mathbb{R}$
\[\label{eq:p}
\min_{x\in\mathbb{R}^d} f(x), \tag{$\Diamond$}
\]
where $f$ could be both nonsmooth and nonconvex (``non''-setting for short). We are particularly interested in algorithms with a finite-time complexity for computing approximately stationary points of Problem~\eqref{eq:p}. 
Note that when $f$ is smooth, it is folkloric that computing an $\epsilon$-stationary point (i.e., $\|\nabla f(x)\|\leq \epsilon$) only requires $O(\epsilon^{-2})$ calls, which is dimension-independent and finite, to the gradient oracle with gradient descent \citep{nemirovskij1983problem}.

In the general Lipschitz ``non''-setting, a widely used generalized subdifferential $\partial f(x)$ is due to \citet[Section 2.1]{clarke1990optimization} (see also \Cref{def:subc}), which reduces to the convex subdifferential (resp.\ gradient) if $f$ is convex (resp.\ smooth). Therefore, by mimicking results in the smooth scenario, it is natural to conjecture that we may be able to design algorithms to compute elements in $\{x: \dist(0,\partial f(x)) \leq \epsilon\}$ in finite time with high probability. However, as shown by \citet[Theorem 5]{zhang2020complexity}, that is impossible for any first-order method. Thus, it is curious to ask: What kind of approximate stationarity concept in the ``non''-setting will admit dimension-free finite-time computation?

\citet{davis2019stochastic} gave a nice answer for the class of $\rho$-weakly convex functions\footnote{Recall $f$ is $\rho$-weakly convex if $f(x)+\frac{\rho}{2}\|x\|^2$ is convex. Weak convexity implies Clarke regularity \citep[Proposition 4.5]{vial1983strong}.} by introducing a notion named near-approximate stationarity (NAS, see \Cref{def:nas}), which is closely related to the gradient of the Moreau envelope of $f$. They showed that a subgradient-type method computes an $(\epsilon,\delta)$-NAS point with $O(\rho^4\delta^{-4}+\epsilon^{-4})$ calls to the subgradient oracle. However, many modern ML models are indeed not weakly convex,\footnote{$g(x)=-\max\{x,0\}$ is not Clarke regular (cf.\ \citep[Definition 2.3.4]{clarke1990optimization}) and not $\rho$-weakly convex for any $\rho\in\mathbb{R}$.} e.g., neural networks with ReLU activation functions. Even worse, 
by extending the Lipschitz hardness results in \citep{kornowski2021oracle}, 
\citet{tian2021hardness} demonstrated that, for any finite $T$, there exists a finte $\rho(T)$ such that, for any $0\leq \epsilon,\delta < \frac{1}{2}$ uniformly, computing an $(\epsilon,\delta)$-NAS point for $\rho(T)$-weakly convex functions within $T$ steps is impossible. 

On the other front, 
starting from the seminal work of \citet{goldstein1977optimization},
a notion named Goldstein approximate stationarity (GAS, see \Cref{def:gas}) exhibits favorable algorithmic consequences. The story begins with an approximation of the Clarke subdifferential $\partial_\delta f(x)$ (see \Cref{def:subg}). If we update iteratively with
\[
x_{k+1} \leftarrow x_k - \delta \cdot g_k / \|g_k\|,
\]
where $g_k \coloneqq \arg\min_{g \in \partial_\delta f(x_k)} \|g\|$ is the minimal norm element in $\partial_\delta f(x_k)$, then we can compute an $(\epsilon,\delta)$-GAS point in $O(\epsilon^{-1}\delta^{-1})$ steps. 
The problem is that obtaining $g_k$ for a general Lipschitz function can be computationally expensive (if possible at all) as there is no known approach to compute $\partial_\delta f(x)$. Therefore, a series of works, e.g., \citep{burke2020gradient}, proposed to build a polyhedral approximation of $\partial_\delta f(x_k)$ via random sampling and compute an approximate $g_k$ by solving a QP in every iteration. However, the number of sampling points needed for meaningful approximation of $\partial_\delta f(x_k)\subseteq \mathbb{R}^d$ is lower bounded by the dimension $d$. Thus, a dimension-free finite-time complexity cannot be achieved within the existing gradient sampling scheme. 

Recently, \citet{zhang2020complexity} introduced a novel algorithm that computes $(\epsilon,\delta)$-GAS points for general Lipschitz functions with a dimension-free finite-time complexity $\widetilde{O}(\epsilon^{-3}\delta^{-1})$. They provided a randomized procedure that compute $\tilde{g}_k: \langle \tilde{g}_k, \partial_\delta f(x_k)\rangle \leq \frac{1}{4} \| \tilde{g}_k \|^2$ with high probability within $\widetilde{O}(\epsilon^{-2})$ oracle calls. However, their assumption on the subgradient oracle is stringent  and hard to be implemented in practice. We restate their assumption below.
\paragraph{Oracle in \citep[Assumption 1(a)]{zhang2020complexity}.} Given $x,d$, the oracle $\mathbb{O}(x,d)$ returns $f(x)$ and a Clarke subgradient $g_x$, such that $g_x \in \partial f(x)$ satisfies $\langle g_x, d\rangle = f^\prime (x; d)$.

Indeed, even computing an \emph{arbitrary} element in $\partial f(x)$ for general Lipschitz functions is highly nontrivial \citep{burke2002approximating,nesterov2005lexicographic,khan2013evaluating,kakade2018provably}, let alone the required subgradient needs to satisfy certain linear equation, which was recognized very early on as impractical \cite{wolfe1975method}.
Such considerations motivate the following question (\textbf{Q1}):
\begin{quote}\centering
	Can we compute GAS points in finite time with \\ a provable and practical algorithm?
\end{quote}

It is notable that GAS is a strictly weaker\footnote{Formally, if $x$ is $(\epsilon,\delta)$-GAS, then $x$ is also $(\epsilon,\delta)$-NAS. See \Cref{def:gas,def:nas}.} stationarity notion than NAS even for continuously differentiable functions \citep[Proposition 1]{kornowski2021oracle} and convex functions (\Cref{prop:cvxpoly}).
However, the computability of NAS is much worse than that of GAS, since finite-time algorithms for NAS only exist if the objective function is $\rho$-weakly convex, which rules out many interesting machine learning models. Thus, it is of interest to ask (\textbf{Q2}): 
\begin{quote}\centering
	Can we compute NAS points in finite time for\\ functions beyond $\rho$-weakly convex, practically? 
\end{quote}

\subsection{Prior Arts}
\vspace{-.1cm}
\paragraph{Asymptotic Analysis.}
The asymptotic computability of Clarke stationary points (i.e., $\{x:0\in\partial f(x)\}$) has been well-understood for quite general functions.
With a  differential inclusion perspective, 
\citet{benaim2005stochastic,majewski2018analysis,davis2020stochastic} studied the asymptotic convergence of subgradient-type methods. In particular, \citet{davis2020stochastic} proved the asymptotic convergence to Clarke stationary points for Whitney stratifiable objective functions, which include deep ReLU neural networks as a special case. \citet{daniilidis2020pathological} demonstrated that the vanilla subgradient method may not converge for general Lipschitz functions even in continuous time.
\vspace{-.2cm}
\paragraph{Finite-Time Analysis.}
In contrast to the asymptotic regime, the finite-time complexity in the general ``non''-setting is still developing.
On the positive side, 
\citep{davis2019proximally,davis2019stochastic} showed that for $\rho$-weakly convex functions, $(\epsilon,\delta)$-NAS is computable with  $O(\rho^4\delta^{-4}+\epsilon^{-4})$ oracle calls. On the negative side,
\citet{kornowski2021oracle} showed that computing NAS for Lipschitz functions in dimension-independent finite time is impossible. \citet{tian2021hardness} sharpened the hardness results for NAS to $\rho$-weakly convex with unbounded $\rho$, thus matching the positive results. For GAS, the gradient sampling scheme \cite{burke2005robust,kiwiel2007convergence,kiwiel2010nonderivative,burke2020gradient} promises finite but dimension-dependent complexity. \citep{zhang2020complexity} reported a novel dimension-independent finite-time algorithm with a impractical  subgradient oracle. A recent concurrent work \cite{davis2021gradient} adopted similar strategy as our \Cref{sec:det} with different algorithmic implementation.
 Another line of research is to exploit structure: \citet{duchi2018stochastic,drusvyatskiy2019efficiency,davis2019stochastic,bolte2018first,beck2020convergence}. In these settings, nonsmoothness and nonconvexity are properly separated making finite-time analysis possible.

\vspace{-.2cm}
\subsection{Contributions}
\vspace{-.1cm}
We highlight the main contributions as follows.
\vspace{-.2cm}
\begin{itemize}
	\item For \textbf{Q1}, we report a practical algorithmic scheme to compute GAS points for general Lipschitz functions with finite-time complexity in both deterministic and stochastic settings.
	\item For \textbf{Q2}, we isolate a new function class within which our new algorithmic scheme computes NAS points in finite time.This goes far beyond existing $\rho$-weakly convex results. Besides, we establish a series of theoretical tools to compute parameters in our new function class.
	\item To demonstrate the wide applicability of the new theoretical framework, we show that $\rho$-margin SVM, $1$-layer, and $2$-layer ReLU neural networks, all being Clarke irregular, satisfy our new conditions.
	\end{itemize}
\paragraph{Notations.} The notation used in this paper is mostly standard in variational analysis.
 $\dist(x,S)\coloneqq\inf_{v\in S} \|v-x\|$; $A \oplus B$ denotes the direct sum of $A$ and $B$; $\mathbb{B}_\epsilon(x)\coloneqq\{v:\|v - x\|\leq \epsilon\}; \mathbb{B}\coloneqq \mathbb{B}_1(0)$;  $\conv S$ is the convex hull of set $S$; $\pi_1 A \coloneqq \{x: \exists y, (x,y) \in A \}$; $A^c$ is the complement of set $A$; $\vvec(X)$ is the vectorization of matrix $X$; $\mathbb{S}^{d-1}\coloneqq \{ x\in\mathbb{R}^d: \|x\|=1\}$.
\section{Preliminaries}\label{sec:prel}
In this section, we introduce the necessary background on variational analysis for Lipschitz functions.
To begin, we recall the following definition of Clarke subdifferential \citep[Theorem 9.61]{rockafellar2009variational}.
\begin{definition}[Clarke subdifferential]\label{def:subc}
	\label{def:subd} Given a point $x$, the Clarke subdifferential of Lipschitz $f$ at $x$ is defined by
	\[
	\partial f(x) \coloneqq \conv\big\{s:\exists x^\prime\! \rightarrow\! x, \nabla f(x^\prime) \textnormal{ exists}, \nabla f(x^\prime)\!\rightarrow\! s\big\}.
	\]
\end{definition}
The following $\delta$-approximation of Clarke subdifferential introduced by \citet{goldstein1977optimization} 
has nice theoretical properties and is convenient for algorithmic developments.
\begin{definition}[Goldstein $\delta$-subdifferential]
	\label{def:subg} Given a point $x$ and $\delta \geq 0$, the Goldstein $\delta$-subdifferential of Lipschitz $f$ at $x$ is defined by
	\[
	\partial_\delta f(x) \coloneqq \conv\big\{\textstyle{\bigcup_{y\in\mathbb{B}_\delta (x)}} \partial f(y) \big\}.
	\]
\end{definition}
We record some useful properties of the Clarke subdifferential and its Goldstein approximation here: 
\begin{fact}[cf. \citet{clarke1990optimization,goldstein1977optimization,zhang2020complexity}] For an $L$-Lipschitz continuous $f$ and $\delta > 0$,
\begin{itemize}
	\item $\partial f(x), \partial_\delta f(x)$ are nonempty, convex, compact;
	\item $\partial f(x) = \cap_{\delta > 0} \cup_{y\in \mathbb{B}_\delta (x)} \partial f(y)$;
	\item $\partial f(x) = \cap_{\delta > 0} \partial_\delta f (x)$;
	\item if $f$ is $C^1$ near $x$, then $\partial f(x) = \{\nabla f(x)\}$;
	\item if $f$ is convex, then $\partial f(x)$ is the convex subdifferential.
\end{itemize}
\end{fact}

We are now ready to introduce two important approximate stationarity notions. We refer the reader to \citep{davis19opt} for a nice expository material.
\begin{definition}[Goldstein approximate stationarity, GAS]\label{def:gas}
	Given a locally Lipschitz function $f:\mathbb{R}^d\rightarrow \mathbb{R}$, we say that $x\in\mathbb{R}^d$ is an $(\epsilon,\delta)$-GAS point if
	\[
	\dist \Big(0, \partial_\delta f(x)\Big) \leq \epsilon.
	\]
\end{definition}

\begin{definition}[near-approximate stationarity, NAS]\label{def:nas}
	Given a locally Lipschitz function $f:\mathbb{R}^d\rightarrow \mathbb{R}$, we say that $x\in\mathbb{R}^d$ is an $(\epsilon,\delta)$-NAS point if
	\[
	\dist \Big(0, \textstyle{\bigcup_{y\in\mathbb{B}_\delta (x)}} \partial f(y)\Big) \leq \epsilon.
	\]
\end{definition}

It is easy to see that if $x$ is NAS, then $x$ is also GAS as $\partial_\delta f(x) \supseteq\cup_{y\in\mathbb{B}_\delta (x)} \partial f(y)$. But the converse does not hold in general, even for continuously differentiable functions.
\begin{fact}[{\citet[Proposition 1]{kornowski2021oracle}}]\label{prop:shamir20p1}
For any $\delta > 0$, there exists a continuously differentiable function $f:\mathbb{R}^2\rightarrow \mathbb{R}$, which is $2\pi$-Lipschitz on $\delta\mathbb{B}$, such that $(0,0)$ is $(0,\delta)$-GAS but $\min_{x\in\delta\mathbb{B}} \|\nabla f(x)\| \geq 1$.
\end{fact}

\Cref{prop:shamir20p1} does not hold for $\rho$-weakly convex functions with sufficiently small $\delta$. Thus, it is still unclear whether NAS and GAS are equivalent assuming $\rho$-weak convexity with finite $\rho\geq 0$. We report below a convex polyhedral version (recall that convexity is $0$-weak convexity), which might be of independent interest.
\begin{restatable}[convex polyhedron]{proposition}{propconvexexamp}\label{prop:cvxpoly}
For any $\delta > 0$, there exists a convex function $f:\mathbb{R}^2\rightarrow \mathbb{R}$, which is $2$-Lipschitz with polyhedral $\partial f$, such that $(0,2\delta)$ is $(0,\delta)$-GAS but $ \min_{y\in\mathbb{B}_\delta \left((0,2\delta)\right)} \dist\big(0, \partial f(y)\big) \geq \frac{2}{5}\sqrt{5}$.
\end{restatable}

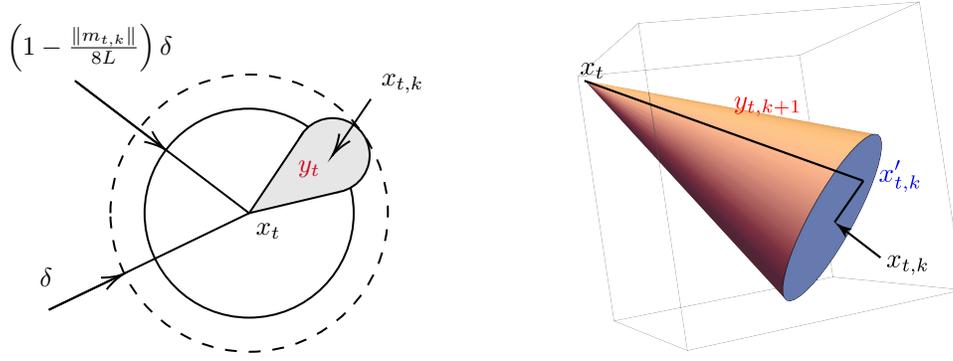
\begin{figure*}[!th]
\centering	
\tikzset{every picture/.style={line width=0.75pt}} %
\begin{tikzpicture}[x=0.75pt,y=0.75pt,yscale=-1,xscale=1]

\draw   (74.29,116.97) .. controls (74.29,87.84) and (97.9,64.23) .. (127.03,64.23) .. controls (156.16,64.23) and (179.77,87.84) .. (179.77,116.97) .. controls (179.77,146.1) and (156.16,169.71) .. (127.03,169.71) .. controls (97.9,169.71) and (74.29,146.1) .. (74.29,116.97) -- cycle ;
\draw  [color={rgb, 255:red, 0; green, 0; blue, 0 }  ,draw opacity=1 ][fill={rgb, 255:red, 230; green, 230; blue, 230 }  ,fill opacity=1 ] (150.3,87.77) .. controls (150.3,77.41) and (158.69,69.02) .. (169.05,69.02) .. controls (179.41,69.02) and (187.8,77.41) .. (187.8,87.77) .. controls (187.8,98.12) and (179.41,106.52) .. (169.05,106.52) .. controls (158.69,106.52) and (150.3,98.12) .. (150.3,87.77) -- cycle ;
\draw  [dash pattern={on 4.5pt off 4.5pt}] (57.06,116.97) .. controls (57.06,78.33) and (88.39,47) .. (127.03,47) .. controls (165.67,47) and (197,78.33) .. (197,116.97) .. controls (197,155.61) and (165.67,186.94) .. (127.03,186.94) .. controls (88.39,186.94) and (57.06,155.61) .. (57.06,116.97) -- cycle ;
\draw [fill={rgb, 255:red, 230; green, 230; blue, 230 }  ,fill opacity=1 ]   (172.8,106.15) -- (127.03,116.97) -- (154.4,76.15) ;
\draw    (39.2,50.15) -- (83.21,83.74) ;
\draw [shift={(84.8,84.95)}, rotate = 217.35] [color={rgb, 255:red, 0; green, 0; blue, 0 }  ][line width=0.75]    (10.93,-3.29) .. controls (6.95,-1.4) and (3.31,-0.3) .. (0,0) .. controls (3.31,0.3) and (6.95,1.4) .. (10.93,3.29)   ;
\draw    (39.2,50.15) -- (127.03,116.97) ;
\draw    (26,165.77) -- (127.03,116.97) ;
\draw    (26,165.77) -- (61.8,148.63) ;
\draw [shift={(63.6,147.77)}, rotate = 154.42] [color={rgb, 255:red, 0; green, 0; blue, 0 }  ][line width=0.75]    (10.93,-3.29) .. controls (6.95,-1.4) and (3.31,-0.3) .. (0,0) .. controls (3.31,0.3) and (6.95,1.4) .. (10.93,3.29)   ;
\draw    (188.4,59.37) -- (170.18,86.12) ;
\draw [shift={(169.05,87.77)}, rotate = 304.27] [color={rgb, 255:red, 0; green, 0; blue, 0 }  ][line width=0.75]    (10.93,-3.29) .. controls (6.95,-1.4) and (3.31,-0.3) .. (0,0) .. controls (3.31,0.3) and (6.95,1.4) .. (10.93,3.29)   ;

\draw (3.2,20) node [anchor=north west][inner sep=0.75pt]  [font=\normalsize]  {$\left( 1-\frac{\| m_{t,k}\| }{8L}\right) \delta $};
\draw (20,143.17) node [anchor=north west][inner sep=0.75pt]  [font=\normalsize]  {$\delta $};
\draw (129.03,120.37) node [anchor=north west][inner sep=0.75pt]  [font=\normalsize]  {$x_{t}$};
\draw (192.2,45.37) node [anchor=north west][inner sep=0.75pt]  [font=\normalsize]  {$x_{t,k}$};
\draw (150.2,89) node [anchor=north west][inner sep=0.75pt]  [font=\normalsize,color={rgb, 255:red, 208; green, 2; blue, 27 }  ,opacity=1 ]  {$y_{t}$};

\end{tikzpicture}
\hspace{5em}
\includegraphics[width=0.3\textwidth]{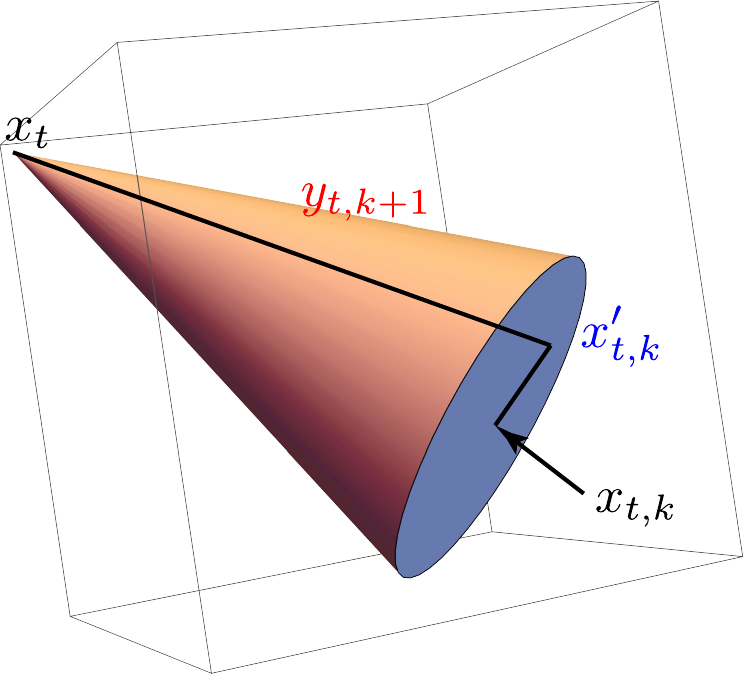}
\caption{Random Conic Perturbation Scheme in \Cref{alg:pingd}.\vspace{-.2cm}}\label{fig:conic}
\end{figure*}

\section{Computing GAS with Practical Oracle by Random Conic Perturbation}
\subsection{Subgradient Oracles}
\begin{assumption}[practical oracle]\label{assu:practial} Given $x$ and Lipschitz continuous $f$:
\begin{itemize}
\item[(a)] In the \textbf{deterministic} setting, if $f$ is differentiable at $x$, then the oracle $\mathbb{O}_\textnormal{d}(x)$ returns a function value $f(x)$ and the gradient $\nabla f(x)$. Otherwise, it sets $\texttt{error}\!=\!1$.
\item [(b)] In the \textbf{stochastic} setting, if $f$ is differentiable at $x$, then the oracle $\mathbb{O}_\textnormal{s}(x)$ returns a stochastic gradient $g_x$ with $\mathbb{E}[g_x \mid \sigma(x)] = \nabla f(x)$ satisfying $\mathbb{E}[\|g_x - \nabla f(x) \|^2 \mid \sigma (x)]\leq \sigma ^2$. Otherwise, it sets $\texttt{error}=1$.	
\end{itemize}

\end{assumption}
Compared with the oracle in \citep[Assumption 1]{zhang2020complexity}, \Cref{assu:practial} only needs to evaluate the gradient $\nabla f(x)$ at differentiable points. Indeed, many modern Automatic Differentiation software (e.g., PyTorch, TensorFlow) can be used as an implementation of \Cref{assu:practial} without worrying about their incorrectness on subgradient evaluation for nonconvex nonsmooth function \citep{kakade2018provably}.
\subsection{Deterministic Setting}\label{sec:det}
In this section, we present a practical algorithm for computing an $(\epsilon,\delta)$-GAS point and establish its finite-time complexity in the deterministic setting. The new algorithm replaces the stringent oracle assumption  in \citep[Assumption 1(a)]{zhang2020complexity} with \Cref{assu:practial}(a). 
\subsubsection{Algorithm}

The main idea is to make use of the almost everywhere differentiability of Lipschitz functions as guaranteed by Rademacher's Theorem. 
 By adopting a random conic perturbation to the uniform sampling direction in \citep[Algorithm 1]{zhang2020complexity}, we have the following \Cref{alg:pingd}, in which the main differences with \citep[Algorithm 1]{zhang2020complexity} are marked in \textcolor{blue}{blue}. See also \Cref{fig:conic}.
\begin{algorithm}[H]
   \caption{Perturbed INGD}
   \label{alg:pingd}
   	\renewcommand{\algorithmicrequire}{\textbf{Input:}}
   \renewcommand{\algorithmicensure}{\textbf{Initialize:}}
\begin{algorithmic}[1]
   \REQUIRE $x_1 \in \mathbb{R}^d$.
   \STATE \textcolor{blue}{Set $\texttt{error}=0$}.
   \FOR{$t\in[T]$}
   \WHILE{$\|m_{t,K}\|>\epsilon$}
   \STATE \textcolor{blue}{Sample $y_{t,1}$ uniformly from $\mathbb{B}_\delta(x_t)$.}
   \STATE Call oracle $\sim, m_{t,1} = \textcolor{blue}{\mathbb{O}_\text{d}(y_{t,1})}$.   
   \FOR{$k\in[K]$}
   \STATE $x_{t,k}=x_t - \textcolor{blue}{\left(1-\frac{\|m_{t,k}\|}{8L}\right)}\cdot\delta \frac{m_{t,k}}{\|m_{t,k}\|}$.
   \IF{$\|m_{t,k}\|\leq \epsilon$}
   \STATE Terminate the algorithm and return $x_t$.
   \ELSIF{$f(x_{t,k})-f(x_t)<-\frac{\delta}{4}\|m_{t,k}\|$}
   \STATE Set $x_{t+1}=x_{t,k}$ and $t=t+1$.
   \STATE Break while-loop.
   \ELSE
   \STATE \textcolor{blue}{Sample $u_{t,k+1}\in\mathbb{R}^{d+1}$ uniformly from $\mathbb{S}^{d}$.}
   \STATE \textcolor{blue}{Let $v_{t,k+1}$ be the first $d$ coordinates of $u_{t,k+1}$.}
   \STATE \textcolor{blue}{$b_{t,k+1}=v_{t,k+1} - \frac{v_{t,k+1}^\top (x_t - x_{t,k})}{\|x_t - x_{t,k}\|^2}\cdot (x_t - x_{t,k})$.}
   \STATE Sample $y_{t,k+1}$ uniformly from $[x_t, \textcolor{blue}{x_{t,k}^\prime}]$, \textcolor{blue}{where $x_{t,k}^\prime\coloneqq x_{t,k}\textcolor{blue}{+\frac{\delta\|m_{t,k}\|}{8L}\cdot b_{t,k+1}}$}.
   \STATE Call oracle $\sim, g_{t,k+1}= \textcolor{blue}{\mathbb{O}_\text{d}(y_{t,k+1})}$.
   \STATE Update $m_{t,k+1}=\beta_{t,k} m_{t,k}+(1-\beta_{t,k})g_{t,k+1}$ 
   \STATE with $\beta_{t,k}=\textcolor{blue}{\frac{8L^3-L^2\|m_{t,k}\|-4L\|m_{t,k}\|^2}{8L^3-L^2\|m_{t,k}\|-\|m_{t,k}\|^3}}$.
   \ENDIF
   \ENDFOR
   \ENDWHILE
   \ENDFOR
   \end{algorithmic}
\end{algorithm}

\subsubsection{Finite-Time Analysis}
The main technical contributions in the analysis are summarized in the following two lemmas.
\begin{restatable}{lemma}{zeroprob} \label{lem:zeroprob}
Let $D\coloneqq\{x: f\textnormal{ is differentiable at }x\}$. Given locally Lipschitz continuous $f$, we have \[
\vspace{-.5em}
\mathbb{P}\Big(\exists (t,k)\in[T]\times [K]: y_{t,k} \in D^c \Big) = 0.
\]

\end{restatable}

\begin{restatable}{lemma}{lemm} \label{lem:m}
	Let $K=\frac{80L^2}{\epsilon^2}$. Given $t\in[T]$, it holds
	\[
	\mathbb{E}\Big[\|m_{t,K}\|^2\Big]\leq \frac{\epsilon^2}{16},
	\]
	where $m_{t,k}=0$ for all $k > k_0$ if the $k$-loop breaks at $(t,k_0)$. Consequently, for any $0\leq\gamma<1$, with probability $1-\gamma$, there are at most $\log(\gamma^{-1})$ restarts of the while loop in the $t$-th iteration.
\end{restatable}

We have the following finite-time guarantee for \Cref{alg:pingd}.
\begin{restatable}{theorem}{thmconicdet}
Let $f$ be $L$-Lipschitz continuous. Then, \Cref{alg:pingd} with $K=\frac{80L^2}{\epsilon^2}$ and $T=\frac{4\Delta}{\epsilon\delta}$ finds an $(\epsilon,\delta)$-GAS point with probability $1-\gamma$ using at most
\[
\frac{320\Delta L^2}{\epsilon^3\delta}\log\left(\frac{4\Delta}{\gamma\epsilon\delta}\right) \qquad\textnormal{oracle calls}
\]
with $\mathbb{P}(\texttt{error}=1)=0$, where $f(x_0) - \inf_x f(x) \leq \Delta$.
\end{restatable}

\subsection{Stochastic Setting}\label{sec:sto}
In this section, we consider the stochastic setting. The new algorithm replaces the stringent oracle assumption  in \citep[Assumption 1(b)]{zhang2020complexity} with \Cref{assu:practial}(b).

\subsubsection{Algorithm}
Technically speaking, the main difference from \citep[Algorithm 2]{zhang2020complexity} lies in the additional perturbation step. We need to carefully choose $\zeta$ to ensure that the iterates are within a $\delta$-ball of some reference point without hurting the convergence. Since $m_t$ is a weighted average of all the stochastic gradients, we need to show that it approximately belongs to the Goldstein $\delta$-subdifferential $\partial_\delta f(x)$ of some reference point $x$.

\textbf{The subtlety when $\norm{m_t} = 0$:} Unlike in the deterministic setting where we can terminate the algorithm if $\norm{m_{t,k}}$ is small, in the stochastic case, $m_t$ is a convex combination of stochastic gradients and thus it does not suffice to terminate the algorithm even if $\norm{m_t} = 0$. The quantity that we aim to minimize is its expectation $\norm{\E{m_t}}\leq \E{\norm{m_t}}$. Due to this subtlety, we cannot let the perturbation size $\zeta_t$ adapt to $\norm{m_t}$ as in the deterministic case: If $\zeta_t = \frac{\omega_1\norm{m_t}}{p\norm{m_t} + \omega_2}$ in Algorithm~\ref{alg:SINGD}, then when $\norm{m_t} = 0$, we have  $y_{t+1} = x_{t+1} = x_t$, and we cannot ensure that $f$ is differentiable at $x_t$ almost surely. We choose a constant $\zeta_t \equiv \zeta$ in Algorithm~\ref{alg:SINGD} instead. In this case, when $\norm{m_t} = 0$, $y_{t+1}$ is sampled from a ball centered at $x_t$.

 By adopting a random conic perturbation to  \citep[Algorithm 2]{zhang2020complexity}, we have the following \Cref{alg:SINGD}, in which the main differences with \citep[Algorithm 2]{zhang2020complexity} are marked in \textcolor{blue}{blue}. 
\begin{algorithm}[H]
	\caption{Perturbed Stochastic INGD}
	\label{alg:SINGD}
	\renewcommand{\algorithmicrequire}{\textbf{Input:}}
	\renewcommand{\algorithmicensure}{\textbf{Initialize:}}
	\begin{algorithmic}[1]
		\REQUIRE $x_1 \in \R^d$.
		\ENSURE $m_1 = \sg_1 = \textcolor{blue}{\so{x_1}}$. Set $\beta = 1 - \frac{\epsilon^2}{64G^2}$, $\textcolor{blue}{K = \frac{1}{\ln{\frac{1}{\beta}}}\ln{\frac{16G}{\epsilon}}}$, $\textcolor{blue}{\omega = \left(\frac{1}{1-\beta} - \frac{1}{\ln{\frac{1}{\beta}}} \right) \ln{\frac{16G}{\epsilon}}}$, $p = \frac{64G^2}{\delta\epsilon^2}\ln{\frac{16G}{\epsilon}}$, $q= \frac{256G^3}{\delta\epsilon^2}\ln{\frac{16G}{\epsilon}}$,  $T = \frac{2^{16}G^3\Delta\ln{\frac{16G}{\epsilon}}}{\epsilon^4\delta}\max\{1, \frac{G\delta}{8\Delta}\} $. 
		\STATE \textcolor{blue}{Set $\texttt{error} = 0$.}
		\FOR{$t \in [T]$}
		\STATE $x_{t+1} = x_t - \eta_t m_t$, where $\eta_t = \frac{1}{p\norm{m_t} + q}$.
		\STATE \textcolor{blue}{Sample $u_{t+1}\in \R^{d+1}$ uniformly from $\mathbb{S}^d$.} 
		\STATE \textcolor{blue}{Let $v_{t+1}\in \R^d$ be the first $d$ coordinates of $u_{t+1}$.}
		\STATE \textcolor{blue}{If $\norm{m_t}>0, b_{t+1} = v_{t+1} - \frac{\innr{v_{t+1}, x_t - x_{t+1}}}{\norm{x_t - x_{t+1}}^2} \cdot (x_t - x_{t+1})$; otherwise,  $b_{t+1} = v_{t+1}$.}
		\STATE Sample $y_{t+1}$ uniformly from $[x_t, x_{t+1} \textcolor{blue}{+ \zeta b_{t+1}}]$, \textcolor{blue}{{where $\zeta = \min\{\frac{\omega}{p}, \frac{\epsilon^2}{510q(L+G)}\}$.}}
		\STATE Call oracle $\sg_{t+1} = \textcolor{blue}{\so{y_{t+1}}}$.
		\STATE $m_{t+1} = \beta m_t + (1 - \beta) \sg_{t+1}$.
		\ENDFOR
		\renewcommand{\algorithmicensure}{\textbf{Output:}}
		\ENSURE $x_{\text{out}} \coloneqq x_{\max\{1, i - K\}}$, where $i \sim \text{Unif}([T])$.
	\end{algorithmic}
	
\end{algorithm}

\vspace{-1em}
\subsubsection{Finite-Time Analysis}
We have the following finite-time guarantee for \Cref{alg:SINGD}, which is similar to \citep[Theorem 10]{zhang2020complexity} but  replaces the stringent oracle assumption  in \citep[Assumption 1(b)]{zhang2020complexity} with \Cref{assu:practial}(b).
\begin{restatable}{theorem}{thmSINGD}
	\label{thm:SINGD} Under \Cref{assu:practial}(b), with probability at least $\frac{3}{5}$, the output of Algorithm \ref{alg:SINGD} satisfies
	$
	\dist(0, \partial_\delta f(x_{\textnormal{out}})) \leq \epsilon
	$
	after at most 
	\[
	\widetilde{O}\left(\frac{G^3\Delta}{\epsilon^4\delta}\right) \qquad \textnormal{oracle calls}
	\]
	with $\mathbb{P}(\texttt{error}=1) = 0$, where $f(x_0) - \inf_x f(x) \leq \Delta$. 
\end{restatable}
\section{Computing NAS by GAS}\label{sec:olc}

In this section, we isolate a new function class within which the new algorithmic scheme can compute near-approximately stationary points in finite time. The new class goes far beyond that of $\rho$-weakly convex functions. We will first introduce the general results, and then several useful calculus rules. In \Cref{sec:application}, we will discuss applications of the new techniques to modern machine learning models.

\subsection{General Results}\label{sec:olcmain}
The main strategy is to compute NAS by GAS.
To this end, we need certain continuity of set-valued subdifferential mapping $\partial f:\mathbb{R}^d\rightrightarrows\mathbb{R}^d$, which should be stronger than upper semicontinuity. A classic notion in set-valued analysis named outer Lipschitz continuity is defined as follows.
\begin{definition}[{\citet[3D]{dontchev2009implicit}}]
	A set-valued mapping $G:\mathbb{R}^d\rightrightarrows\mathbb{R}^d$ is outer Lipschitz continuous (OLC) at $\bar{y}$ relative to a set $D$ if $\bar{y}\in D \subset \dom G$, $G(\bar{y})$ is a closed set, and there is a constant $\kappa\geq 0$ along with a neighborhood $V$ of $\bar{y}$ such that
	\[
	G(y) \subseteq G(\bar{y}) + \kappa \|y - \bar{y}\|\mathbb{B}, \qquad \forall y \in V\cap D.
	\]
\end{definition}

OLC is a weaker notion than Lipschitz continuity even for a single-valued mapping $G:\mathbb{R}\rightarrow\mathbb{R}$. See \citep[Example 2.4(a)]{lewis2010lipschitz}.
	However, for our purposes, OLC is not sufficient since by the classic result of \citet{robinson1981some} the bad function in \Cref{prop:cvxpoly} is OLC.

The following modified OLC notion for set-valued mapping is new and central in our development, which allows us to have a Lipschitz-type control of $G:\mathbb{R}^d\rightrightarrows\mathbb{R}^d$ from above within a constant-size neighborhood (see also \Cref{fig:polc}).
\begin{definition}[$(\delta,\eta,\kappa)$-outer Lipschitz continuous]\label{def:polc}
	A set-valued mapping $G:\mathbb{R}^d\rightrightarrows\mathbb{R}^d$ is $(\delta,\eta,\kappa)$-OLC on $S$ if for any $x\in S$, there exists a \textbf{pivot} $y\in \mathbb{B}_\delta(x)\cap S$ such that $G$ is $\kappa$-OLC on $\mathbb{B}_\eta(x)\cap S$. In other words, for all $ x \in S$, there exists a $y \in \mathbb{B}_\delta(x)\cap S$ such that
\[
  G(z) \subseteq G(y) + \kappa \|y-z\| \mathbb{B},\ \forall z \in \mathbb{B}_\eta(x)\cap S.
\]
Besides, we call $P^G: x \rightarrow y$ the \textbf{pivot mapping} of $G$.
\end{definition}

\begin{figure}[h]
\centering

\tikzset{every picture/.style={line width=1pt}} %

\begin{tikzpicture}[x=0.75pt,y=0.75pt,yscale=-1,xscale=1]

\draw  [fill={rgb, 255:red, 240; green, 240; blue, 240 }  ,fill opacity=1 ] (91.3,152) .. controls (91.3,133.39) and (106.39,118.3) .. (125,118.3) .. controls (143.61,118.3) and (158.7,133.39) .. (158.7,152) .. controls (158.7,170.61) and (143.61,185.7) .. (125,185.7) .. controls (106.39,185.7) and (91.3,170.61) .. (91.3,152) -- cycle ;
\draw  [dash pattern={on 4.5pt off 4.5pt}] (73.13,152) .. controls (73.13,123.35) and (96.35,100.13) .. (125,100.13) .. controls (153.65,100.13) and (176.88,123.35) .. (176.88,152) .. controls (176.88,180.65) and (153.65,203.88) .. (125,203.88) .. controls (96.35,203.88) and (73.13,180.65) .. (73.13,152) -- cycle ;
\draw    (154.2,85) -- (138.62,119.58) ;
\draw [shift={(137.8,121.4)}, rotate = 294.25] [color={rgb, 255:red, 0; green, 0; blue, 0 }  ][line width=0.75]    (10.93,-3.29) .. controls (6.95,-1.4) and (3.31,-0.3) .. (0,0) .. controls (3.31,0.3) and (6.95,1.4) .. (10.93,3.29)   ;
\draw    (53.4,204.2) -- (80.98,184.18) ;
\draw [shift={(82.6,183)}, rotate = 144.02] [color={rgb, 255:red, 0; green, 0; blue, 0 }  ][line width=0.75]    (10.93,-3.29) .. controls (6.95,-1.4) and (3.31,-0.3) .. (0,0) .. controls (3.31,0.3) and (6.95,1.4) .. (10.93,3.29)   ;
\draw  [fill={rgb, 255:red, 0; green, 0; blue, 0 }  ,fill opacity=1 ] (165.2,182.3) .. controls (165.2,181.14) and (166.14,180.2) .. (167.3,180.2) .. controls (168.46,180.2) and (169.4,181.14) .. (169.4,182.3) .. controls (169.4,183.46) and (168.46,184.4) .. (167.3,184.4) .. controls (166.14,184.4) and (165.2,183.46) .. (165.2,182.3) -- cycle ;

\draw (120,148) node [anchor=north west][inner sep=0.75pt]    {$x$};
\draw (157.6,72.6) node [anchor=north west][inner sep=0.75pt]    {$\eta $};
\draw (43.2,182) node [anchor=north west][inner sep=0.75pt]    {$\delta $};
\draw (171.6,183.6) node [anchor=north west][inner sep=0.75pt]    {$y=P( x)$};

\end{tikzpicture}
\caption{$(\delta,\eta,\kappa)$-Outer Lipschitz Continuity in \Cref{def:polc}.}\label{fig:polc}
\end{figure}
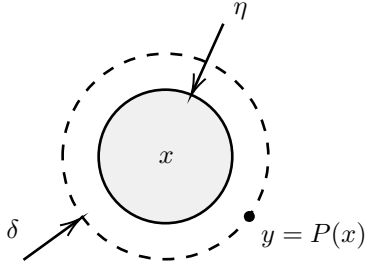
\begin{remark}
A natural question about \Cref{def:polc} is why we want to set $\eta$ and $\delta$ to different values. In other words, why $(\delta,\delta,\kappa)$-OLC is not sufficient. 
Consider the convex function $g(x,y)\coloneqq \max\{2x,-2x,y\}$, which is the bad function in the proof of \Cref{prop:cvxpoly}. It is easy to see that $\partial g$ is polyhedral and OLC by \citep{robinson1981some}.
However, for any $\delta > 0$, $\partial g$ is not $(\delta,\delta,\kappa)$-OLC at $(0,2\delta), \forall \kappa \geq 0$. Thus, even for an OLC mapping $\partial g$, we cannot promise $\exists \delta > 0$ such that $\partial g$ is $(\delta,\delta,\kappa)$-OLC at certain $x\in\dom \partial g$. Instead, we will show in \Cref{thm:generic} that if $\partial g$ is OLC and $S$ is compact, then $\forall \delta>0, \exists \eta > 0$ such that $\partial g$ is $(\delta,\eta,\kappa)$-OLC on $S$.
\end{remark}

We are now ready for the main theorem of this section:
\begin{restatable}[NAS by GAS]{theorem}{thmmain}\label{thm:eps-del-main}
	For a Lipschitz continuous $f$, suppose that $\partial f:\mathbb{R}^d\rightrightarrows\mathbb{R}^d$ is $(\delta,\eta,\kappa)$-OLC. If $x$ is $(\epsilon,\eta)$-GAS, then $x$ is $\big(\epsilon + \kappa (\delta+\eta),\delta\big)$-NAS.
\end{restatable}

It is natural to ask what function class admits a $(\delta,\eta,\kappa)$-OLC subdifferential.
\begin{restatable}{theorem}{thmgen}\label{thm:generic}
	Let $\delta > 0$ and $\partial f:\mathbb{R}^d\rightrightarrows\mathbb{R}^d$  be  $\kappa$-OLC. For any compact set $S$, there exists an $\eta \in (0,\delta]$ such that $\partial f$ is $(\delta,\eta,\kappa)$-OLC on $S$.
\end{restatable}

\begin{remark}
	If the set of $(\epsilon,\delta)$-GAS points is bounded and $\partial f$ is OLC, then we can use \Cref{thm:generic} and \Cref{thm:eps-del-main} to guarantee NAS from GAS. Note that functions with OLC subdifferential have been widely studied in the variational analysis literature. For example, 
	$\partial f$ with a finite union of convex polyhedral graph \citep{robinson1981some} is OLC. If $f$ is an $\ell c$-stable function \citep[Theorem 2]{bednavrik2013lipschitz}, then $\partial f$ is OLC.
	\end{remark}

	Given an OLC mapping $\partial f$ and a constant $\delta>0$, it is generally hard to estimate the constant $\eta$ as guaranteed by Theorem 4.5. However, its value is needed for the stopping rules of \Cref{alg:pingd}. In the next subsection, we provide several useful calculus rules to compute the parameter $\eta$ explicitly.

\subsection{Calculus of $(\delta,\eta,\kappa)$-Outer Lipschitz Continuity}\label{sec:rules}
In this section, we establish a series of calculus rules to verify and compute the parameters in \Cref{def:polc}. We first introduce four rules that have taken the subdifferential calculus rules\footnote{We note here that the validity of subdifferential chain rules for nonconvex nonsmooth functions is highly non-trivial. See, for example, \citep[Chapter 10]{rockafellar2009variational}.} of $f$ into consideration.

\begin{restatable}[smooth regularization]{proposition}{propcaclsmooth}\label{prop:calc-smooth}
	Suppose that $f:\mathbb{R}^d\rightarrow \mathbb{R}$ has a $(\delta,\eta,\kappa)$-OLC $\partial f$ and $g:\mathbb{R}^d\rightarrow \mathbb{R}$ is differentiable with a $\beta$-Lipschitz gradient $\nabla g$. Then $\partial (f+g)$ is  $(\delta,\eta,\beta+\kappa)$-OLC.
\end{restatable}
\begin{restatable}[separable sum]{proposition}{propcalcsum}\label{prop:calc-sum}
	Suppose, for any $i\in[m]$, that $f_i:\mathbb{R}^{d_i}\rightarrow \mathbb{R}$ has a $(\delta_i,\eta_i,\kappa_i)$-OLC $\partial f_i$. Let $f(x)\coloneqq\sum_{i=1}^m f_i(x_i)$, where $x\coloneqq \bigoplus_{i=1}^m x_i$. Then, $\partial f$ is  $(\delta,\eta,\kappa)$-OLC with
	\[
	\delta = \sqrt{\sum_{i=1}^m \delta_i^2},\quad \eta = \min_{i \in [m]} \eta_i, \quad \kappa = \sqrt{\sum_{i=1}^m \kappa_i^2}.
	\]
\end{restatable}

\begin{restatable}[linear composition]{proposition}{propcalclinear}\label{prop:calc-linear}
	Suppose that $f:\mathbb{R}^n\rightarrow \mathbb{R}$ has a $(\delta,\eta,\kappa)$-OLC $\partial f$ and $A\in\mathbb{R}^{n\times d}$ is surjective. Then, $\partial (f\circ A)$ is  $\left(\delta\|A^\dag\|,\frac{\eta}{\|A\|},\kappa \|A\|^2\right)$-outer Lipschitz continuous.
\end{restatable}
\begin{restatable}[rescaling]{proposition}{propcalcrescal}\label{prop:calc-rescaling}
	Suppose that the $L_1$-Lipschitz $f:\mathbb{R}^n\rightarrow \mathbb{R}$ has a $(\delta,\eta,\kappa)$-OLC $\partial f$ and $g:\mathbb{R}\rightarrow \mathbb{R}$ is $L_2$-Lipschitz and $\beta$-smooth. Then, $\partial (g\circ f)$ is  $\left(\delta,\eta,\beta L_1 + \kappa L_2\right)$-OLC.
\end{restatable}

Then, we introduce a partial sum rule, which is powerful but needs to be used in conjunction with certain subdifferential calculus rules (e.g., assuming Clarke regularity). The following rules are crucial in the 2-layer ReLU neural network example (see \Cref{sec:2relu}).

\begin{restatable}[sum]{proposition}{propcalcsumshared}\label{prop:calc-psum}
	Suppose, for any $ i\in[m]$, that $G_i:\mathbb{R}^{d}\rightrightarrows \mathbb{R}^d$ is $(\delta_i,\eta_i,\kappa_i)$-OLC with a shared pivot mapping $P:\mathbb{R}^d\rightarrow\mathbb{R}^d$. Let $G(x)\coloneqq\sum_{i=1}^m G_i(x)$. Then, $G$ is  $(\delta,\eta,\kappa)$-OLC with
	\[
	\delta = \min_{i \in [m]} \delta_i,\quad \eta = \min_{i \in [m]} \eta_i, \quad \kappa = \sum_{i=1}^m \kappa_i.
	\]
\end{restatable}

\begin{restatable}[partially separable sum]{corollary}{propcalcpartial}\label{prop:calc-partial}
	Suppose, for any $ i\in[m]$, that $G_i:\mathbb{R}^{d_0}\times\mathbb{R}^{d_i}\rightrightarrows \mathbb{R}^{d_0}\times\mathbb{R}^{d_i}$ is $(\delta_i,\eta_i,\kappa_i)$-OLC with a partially shared pivot mapping $P_i:\mathbb{R}^{d_0}\times\mathbb{R}^{d_i}\rightarrow \mathbb{R}^{d_0}\times\mathbb{R}^{d_i}$, such that $\pi_1 \circ P_i(x_0, x_i) = \pi_1 \circ P_1(x_0, x_1), \forall i \in [m]$. Let $G(x)\coloneqq\sum_{i=1}^m G_i(x_0, x_i)$, where $x\coloneqq \bigoplus_{i=0}^m x_i$. Then, $G$ is  $(\delta,\eta,\kappa)$-OLC with
	\[
	\delta = \sqrt{\sum_{i=1}^m \delta_i^2},\quad \eta = \min_{i \in [m]} \eta_i, \quad \kappa = \sum_{i=1}^m \kappa_i.
	\]
\end{restatable}

\subsection{Discussion}
We record here a recipe to prove $(\delta,\eta,\kappa)$-OLC from scratch, which when combined with the calculus rules in this section forms a toolbox for determining the parameters $(\delta,\eta,\kappa)$.
\begin{itemize}
	\item[S1.] Construct pivot mapping $P:\mathbb{R}^d\rightarrow\mathbb{R}^d$.
	\item[S2.] Verify $\|x - P(x)\| \leq \delta$ for all $x\in\mathbb{R}^d$.
	\item[S3.] Prove that for all $x\in\mathbb{R}^d$, it holds
	\[
	G(z) \subseteq G(P(x)) + \kappa \|z - P(x)\| \mathbb{B},\ \forall z \in \mathbb{B}_\eta(x)\cap S.
	\]
\end{itemize}
We will provide concrete examples in \Cref{sec:application}.

\section{Applications}\label{sec:application}
To demonstrate the wide applicability of the new theoretical framework, we discuss examples in machine learning, namely $\rho$-margin SVM, $1$-layer, and $2$-layer ReLU NN, all being Clarke irregular and not weakly convex. We show that all these examples are subdifferential $(\delta,\eta,\kappa)$-OLC, where the parameters $(\delta,\eta,\kappa)$ can be determined via  the calculus rules in \Cref{sec:rules}.
\subsection{$\rho$-Margin loss SVM}
We aim to solve
\[\label{eq:msvm}
\min_{w\in\mathbb{R}^d} F(w)\coloneqq \frac{1}{2} \|w\|^2 + \sum_{i=1}^n \phi_\rho (z_i^\top w), \tag{$\rho$-MSVM}
\]
where $ \phi_{\rho}(u) \coloneqq \min\left(1,\max\left(0, 1-\frac{u}{\rho}\right)\right)$.

The goal is to compute $(\epsilon,\delta)$-NAS points for Problem~\eqref{eq:msvm} by computing $(\epsilon^\prime,\delta^\prime)$-GAS points.
We note that the $\rho$-Margin loss SVM in Problem~\eqref{eq:msvm}  and its $\rho=1$ version, also known as ramp loss SVM, have been widely recognized in the operations research \cite{brooks2011support,carrizosa2014heuristic,wang2021proximal,tian2022computing}, statistics \cite{shen2003psi,wu2007robust,liu2005multicategory}, and machine learning \cite{huang2014ramp,david2011generalization,collobert2006large,collobert2006trading,ertekin2010nonconvex,suzumura2017homotopy,maibing2015computational} communities as providing better robustness against data outliers than the vanilla SVM. The general $\rho$-version can be found in the learning theory textbook \citep[Corollary 5.11]{mohri2018foundations}.

It is elementary to see that $\partial \phi_\rho$ is $(\delta,\delta,0)$-OLC for any $0<\delta\leq \frac{\rho}{2}$ with pivot mapping $P^{\phi_\rho}:\mathbb{R}\rightarrow\mathbb{R}$ defined by
\[
P^{\phi_\rho}(x) \coloneqq 
\left\{ \begin{array}{rcl}
	0 & \mbox{for} 
	& |x|\leq \frac{\rho}{2}, \\ 
	\rho & \mbox{for} & |x-\rho| < \frac{\rho}{2}, \\
	x & \mbox{for} & \text{otherwise}.
\end{array}\right.
\]

Let $\Phi_\rho(y)\coloneqq\sum_{i=1}^n \phi_\rho(y_i)$. Then, by \Cref{prop:calc-sum}, $\partial \Phi_\rho$ is $(\sqrt{n}\delta,\delta,0)$-OLC. Assuming that $Z\in\mathbb{R}^{n\times d}$ is surjective, by \Cref{prop:calc-linear}, $\partial (\Phi_\rho \circ Z)$ is $\left(\sqrt{n}\delta\|Z^\dag\|,\frac{\delta}{\|Z\|},0\right)$-OLC. Using \Cref{prop:calc-smooth}, $\partial F$ is $\left(\sqrt{n}\delta\|Z^\dag\|,\frac{\delta}{\|Z\|},1\right)$-OLC. By \Cref{thm:eps-del-main}, if $x$ is  $\left(\epsilon, \frac{\delta}{\|Z\|}\right)$-GAS, then it is also $\left(\epsilon + \left(\sqrt{n}\|Z^\dag\|+\frac{1}{\|Z\|}\right)\delta, \sqrt{n}\|Z^\dag\|\delta\right)$-NAS.
Let the condition number of $Z$ be $\kappa(K)\coloneqq \|Z^\dag\|\|Z\|$. In other words, to compute an $(\epsilon, \delta)$-NAS point, it is sufficient to have an $(\epsilon^\prime, \delta^\prime)$-GAS point, where (in a dimension-free manner)
\[
\epsilon^\prime \leq \frac{\epsilon}{2} \textnormal{ and }  \delta^{\prime} \leq \min\left\{\frac{\delta}{\sqrt{n}\kappa(Z)}, \frac{\epsilon}{2\sqrt{n}\kappa(Z)+2},\frac{\rho}{2\|Z\|} \right\}.
\]

\subsection{Shallow ReLU Neural Network.}
In this subsection, we will discuss the computation of $(\epsilon,\delta)$-NAS points for shallow ReLU neural networks.
For simplicity, we will not trace explicitly the constants $(\delta,\eta,\kappa)$ in this subsection. Instead, we will say that $f$ is \emph{subdifferentially OLC trackable} if the parameters $(\delta,\eta,\kappa)$ of $\partial f$ can be determined by the calculus rules in \Cref{sec:rules}.

Recently, finite-time convergence of neural networks in the overparameterized regime has been extensively studied 
\citep{jacot2018neural,chizat2018lazy,du2018gradient,arora2019fine,du2019gradient,zou2020gradient}. For the underparameterized regime, the asymptotic convergence of ReLU neural network is analyzed in the continuous-time gradient flow sense \citep{eberle2021existence,jentzen2021existence}. However, it is still unclear what convergence guarantee we can have for (potentially underparameterized) ReLU neural networks within finite time as they are not weakly convex and the finite-time analyses in \cite{davis2019proximally,davis2019stochastic} are inapplicable.

\subsubsection{$1$-Layer ReLU Neural Network}
We first investigate the easy case, that is, the $2$-layer ReLU neural network with the weights of the second layer fixed. It is notable that we will not impose any assumption on the number of hidden nodes $m$.

Let $\sigma(u)\coloneqq \max\{u,0\}$. Setting pivot $y\coloneqq \mathbb{1}_{|u|>\delta} u$, it is elementary to see that $\partial \sigma:\mathbb{R}\rightrightarrows\mathbb{R}$ is $(\delta,\delta,0)$-OLC for any $\delta > 0$. Similarly, $\partial (-\sigma)$ is $(\delta,\delta,0)$-OLC for any $\delta > 0$.
We aim to solve
\[
\min_{W \in \mathbb{R}^{d\times m}}\! F(W) \coloneqq\! \sum_{i=1}^n \ell\!\left(\!y_i, \sum_{j=1}^m (-1)^j\sigma \left(w_j^\top x_i\right)\!\right) + R(W).
\]
Suppose that the regularization term $R:\mathbb{R}^{d\times m} \rightarrow \mathbb{R}$ is smooth. Let $h:\mathbb{R}^{m}\rightarrow\mathbb{R}$ be given by $h(u)\coloneqq\sum_{j=1}^m (-1)^j\sigma (u_j)$.
By \Cref{prop:calc-sum}, $h$ is subdifferentially OLC trackable. Let $\ell_i(u) = \ell(y_i, u)$ and assume that $\ell_i$ is Lipschitz and smooth. Let $f:\mathbb{R}^{mn}\rightarrow \mathbb{R}$ be given by $f(\vvec(U))\coloneqq\sum_{i=1}^n \ell_i \circ h(u_{i})$, where $u_i \in \mathbb{R}^m, \forall i \in [n]$ and $U\in\mathbb{R}^{m\times n}$. With \Cref{prop:calc-rescaling} and using \Cref{prop:calc-sum} again, $f$ is subdifferentially OLC trackable.  We assume that the data $X\in\mathbb{R}^{n\times d}$ is surjective, which holds in many modern high-dimensional machine learning scenarios. Let $x_i\in\mathbb{R}^{1\times d}$ be the $i$-th row of $X$. We define $X_i\in\mathbb{R}^{m\times md}$ and $X_\textnormal{big}\in\mathbb{R}^{mn\times md}$ as
\[
X_i \coloneqq \begin{bmatrix} 
	x_i  & & & \\
		& x_i  & & \\
	& &  \ddots & \\
	&   &   & x_i  
\end{bmatrix},\qquad
X_\textnormal{big} \coloneqq \begin{bmatrix} 
	X_1 \\
	X_2 \\
	\vdots\\
	X_n
\end{bmatrix}.
\]
As $X$ is surjective, $X_\textnormal{big}$ is surjective. Using \Cref{prop:calc-linear}, we have $f(X_\textnormal{big}\vvec(W))$ is subdifferentially OLC trackable, where $\vvec(W) \in \mathbb{R}^{md}$. By $F(W) = f(X_\textnormal{big}\vvec(W))+R(W)$ and  \Cref{prop:calc-smooth}, $F(W)$ is subdifferentially OLC trackable. 

\subsubsection{$2$-Layer ReLU Neural Network}\label{sec:2relu}
Let $\varrho(a,b)\coloneqq a\cdot \max\{b,0\}$. 
We aim to solve
\[
\min_{\substack{W \in \mathbb{R}^{d\times m}\\ a\in\mathbb{R}^m}}\! F(W,a)\! \coloneqq\! \sum_{i=1}^n \ell\left(\!y_i, \sum_{j=1}^m \varrho(a_j,w_j^\top x_i)\!\right)\! + R(W,a)
\]
with surjective $X \in \mathbb{R}^{n\times d}$ and smooth regularization term $R:\mathbb{R}^{d\times m} \times \mathbb{R}^m \rightarrow \mathbb{R}$.

Compared with the $1$-layer case, the main difficulty in the analysis is due to the inseparability of $\{a_j\}_{j\in[m]}$, as one cannot apply the subdifferential chain rule and OLC calculus rules directly. To cope with this, we need the partial separable rule in \Cref{prop:calc-partial} and a partially differentiable sum rule in \Cref{prop:pchain}, which might be of independent interest.
To begin, we have the following subdifferential characterization of $\partial \varrho:\mathbb{R}^2 \rightrightarrows \mathbb{R}^2$:
\begin{restatable}{claim}{clmsubdrho}\label{clm:subd-rho} For $\varrho(u_1,u_2)\coloneqq u_1\cdot \max\{u_2,0\}$, it holds
\[
\partial \varrho (u_1, u_2) = 
\left\{ \begin{array}{rcl}
	(u_2, u_1) & \mbox{for} 
	& u_2 > 0, \\ 
	(0,0) & \mbox{for} & u_2 < 0, \\
	(0, \conv\{0,u_1\}) & \mbox{for} & u_2 = 0.
\end{array}\right.
\]
\end{restatable}

Then, we investigate the continuity of $\partial \varrho$. Given any $\delta>0$, $x\in\mathbb{R}^2$, and $z\in\mathbb{B}_\delta(x)$, we consider the following cases:
\begin{itemize}
	\item If $|x_2| > \delta$, let $y = x$.
	\begin{itemize}
		\item If $y_2>0$, then $z_2 > 0$. We have
		$
		\partial \varrho (z)= (z_2, z_1) \!\subseteq\! (y_2, y_1) \!+\! \|y-z\|\mathbb{B}\! =\! \partial \varrho (y) \!+\! \|y-z\|\mathbb{B}.
		$
		\item If $y_2 < 0$, then $z_2 < 0$. We have $\partial \varrho (z) = (0,0) = \partial \varrho (y)$.
	\end{itemize}
	\item If $0\leq |x_2| \leq \delta$, let $y = (x_1, 0)$. It is easy to see that $\|y-x\| = |x_2| \leq \delta$.
		\begin{itemize}
		\item If $z_2>0$, we have
		$
		\partial \varrho (z)= (z_2, z_1) \subseteq (0, y_1) + \|y-z\|\mathbb{B} \subseteq \partial \sigma (y) + \|y-z\|\mathbb{B}.
		$
		\item If $z_2<0$, we have
		$
		\partial \varrho (z)= (0, 0) \subseteq (0, 0) + \|y-z\|\mathbb{B} \subseteq \partial \sigma (y) + \|y-z\|\mathbb{B}.
		$
		\item If $z_2 = 0$, we have $\partial \varrho (z) = (0,\conv\{0,z_1\}) \subseteq (0, \conv\{0,y_1\}]) + \|y-z\|\mathbb{B} = \partial \sigma (y) + \|y-z\|\mathbb{B}$.
		\end{itemize}
\end{itemize}
Therefore, for any $\delta >0$, $\partial \varrho$ is $(\delta,\delta,1)$-OLC with pivot mapping 
$P^{\partial\varrho}:\mathbb{R}^2\rightarrow\mathbb{R}^2$ defined by
\[
P^{\partial\varrho}\big((x_1, x_2)\big) \coloneqq 
\left\{ \begin{array}{ccl}
	(x_1,0) & \mbox{for} 
	& |x|\leq \delta, \\ 
	(x_1, x_2) &  & \text{otherwise}.
\end{array}\right.
\]
It is easy to see that $\pi_1\circ P^{\partial\varrho}\big((x_1, x_2)\big)$ is independent of $x_2$.
Let $h_i:\mathbb{R}^m \times \mathbb{R}^m \rightarrow \mathbb{R},\forall i \in [n]$ be defined by
\[
h_i(a, u_i)\coloneqq \ell_i \left( \sum_{j=1}^m \varrho(a_j, u_{ij}) \right).
\]
Then, by the choices of pivots in the proof of \Cref{prop:calc-psum,prop:calc-rescaling}, $\partial h_i$ is subdifferentially OLC trackable
with pivot mapping 
$P^{\partial h_i}:\mathbb{R}^m \times \mathbb{R}^m\rightarrow\mathbb{R}^m \times \mathbb{R}^m$ defined by
$
P^{\partial h_i}\big((a, u_i)\big) \coloneqq (a, \tilde{u}_i),
$
where
\[
\tilde{u}_{ij} \coloneqq 
\left\{ \begin{array}{ccl}
	0 & \mbox{for} 
	& |u_{ij}|\leq \delta, \\ 
	u_{ij} & & \text{otherwise}.
\end{array}\right.
\]
Therefore, $\{\partial h_i\}_{i\in [m]}$ partially shares the pivot mapping $P^{\partial h_i}$ on the first argument, i.e., $\pi_1 \circ P^{\partial h_i} \big((x_0, x_i)\big) = \pi_1 \circ P^{\partial h_1} \big((x_0, x_1)\big), \forall i \in [n]$.
Let $f(a, U)\coloneqq \sum_{i=1}^n h_i (a, u_i)$. By \Cref{prop:calc-partial}, $\sum_{i=1}^n \partial h_i$ is subdifferentially OLC trackable. To proceed, we need the following chain rule, whose proof is technical and might be of independent interest.
\begin{restatable}[partially differentiable sum rule]{proposition}{clmchain}\label{prop:pchain} It holds
\[
\partial f(a, U) = \sum_{i=1}^n \partial h_i (a, u_i).
\]
\end{restatable}
Then, $\partial f$ is subdifferentially OLC trackable. Suppose that the data $X\in\mathbb{R}^{n\times d}$ is surjective. Let $x_i$ be the $i$-th row of $X$. We define $\theta\in\mathbb{R}^{m+md}$ and $X_\textnormal{huge}\in\mathbb{R}^{(m+mn)\times (m+md)}$ as
\[
\theta \coloneqq \begin{bmatrix} 
	a   \\
	\vvec(W)
\end{bmatrix},\qquad
X_\textnormal{huge} \coloneqq \begin{bmatrix} 
	I_m  &   \\
	      & X_\textnormal{big}  
\end{bmatrix}.
\]
As $X$ is surjective, $X_\textnormal{huge}$ is surjective. Using \Cref{prop:calc-linear}, we have $f(X_\textnormal{huge}\theta)$ is subdifferentially OLC trackable. By $F(W,a) = f(X_\textnormal{huge}\theta) + R(W,a)$ and  \Cref{prop:calc-smooth}, $F(W,a)$ is subdifferentially OLC trackable.

\section{Closing Remarks}
In this paper, we report a practical algorithmic scheme to compute GAS points for general Lipschitz functions with finite-time complexity.
We also isolate a new function class for which our scheme computes NAS points in finite time. Besides, we establish a series of theoretical tools to compute parameters in our new function class.
To demonstrate the wide applicability of our new theoretical framework, we discuss modern machine learning models and show that they satisfy our new conditions.
We hope that our results can be beneficial to the understanding of finite-time complexity of sharper approximate stationarity for Lipschitz continuous ``non''-problems.
 An intriguing further direction is to apply the new analytical framework to other nonconvex nonsmooth problems. Extending the calculus rules in \Cref{sec:rules} or refining the modified OLC notion in \Cref{def:polc} would also be interesting.

\bibliography{gas-lip}
\bibliographystyle{icml2022}

\newpage
\appendix
\onecolumn
\section{Proofs of \Cref{sec:prel}}

\begin{figure*}[!th]
\centering	
\includegraphics[width=0.35\textwidth]{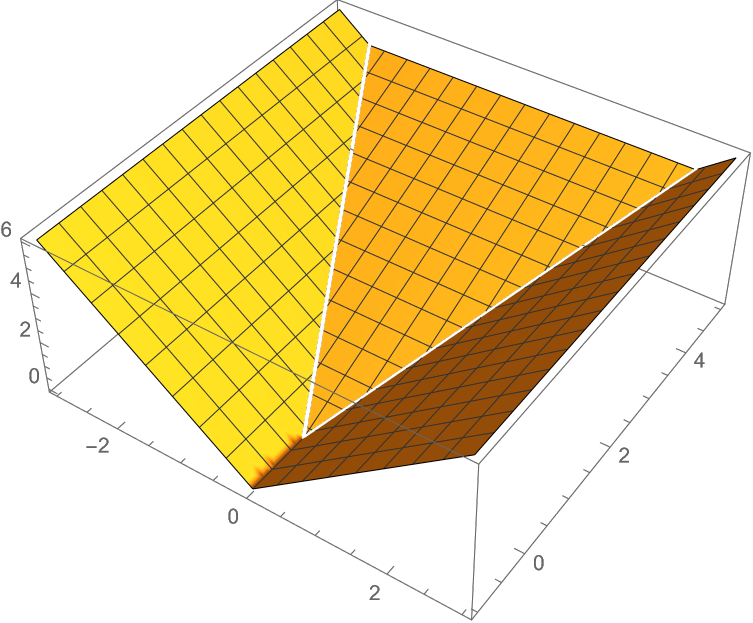}
\caption{The function used in the proof of \Cref{prop:cvxpoly}. }\label{fig:cvxpoly}
\end{figure*}

\propconvexexamp*
\begin{proof}
	Fixing some $\delta > 0$, consider the function (see also \Cref{fig:cvxpoly}), whose convexity is obvious,
	\[
	f(x,y)\coloneqq \max\{2x, -2x, y\}.
	\]
	Note that, by $(\pm \delta,2\delta) \in \mathbb{B}_\delta \big( (0,2\delta) \big)$, it holds
	\[
	(0,0) = \frac{1}{2} (-2,0) + \frac{1}{2} (2,0) \in \frac{1}{2}\partial f\big( (-\delta,2\delta) \big) + \frac{1}{2}\partial f\big( (\delta,2\delta) \big) \subseteq \partial_\delta f\big( (0,2\delta) ).
	\]
	Besides, as $(0,0) \notin \mathbb{B}_\delta \big( (0,2\delta) \big)$, it is elementary to see
	\[
	\dist\Big(0, {\textstyle \bigcup_{y\in\mathbb{B}_\delta ((0,2\delta))}}\partial f (y)\Big) \geq \min_{0\leq \lambda \leq 1} \| (2\lambda,0) + (0, 1-\lambda) \| = \frac{2}{5}\sqrt{5},
	\]
	as required.
\end{proof}

\section{Proofs of \Cref{sec:det}}

\zeroprob*
\begin{proof}
Fix $(t,k)\in[T]\times [K]$. Let 
\[
\begin{aligned}
	S_1 &\coloneqq \left\{ (\lambda, \xi): \lambda \in [0,1], \xi \in \mathbb{R}^{d-1}, \|\xi\| \leq 1 \right\}, \\
	S_2 & \coloneqq \left\{ y \in \mathbb{R}^d : y=x_t + \lambda \left(x_{t,k}-x_t + \frac{\delta \|m_{t,k}\|}{8L} \cdot b_{t,k+1}\right), \lambda \in [0,1], \|b_{t,k+1}\| \leq 1, b_{t,k+1}^\top (x_{t,k}-x_t)=0  \right\}.
\end{aligned}
\]
Let $X^\bot \in \mathbb{R}^{d\times d-1}$ be an orthonormal basis of $\mathop{\textnormal{span}}(x_{t,k}-x_t)^\bot$. We define the following isomorphism:
	\[
	\begin{aligned}
		T: S_1 &\longrightarrow S_2 \\
		(\lambda, \xi) & \longrightarrow y_{t,k+1} \coloneqq x_t + \lambda \left(x_{t,k} - x_t  + \frac{\delta \|m_{t,k}\|}{8L} \cdot X^\bot \xi \right).
	\end{aligned}
	\]
	Then, by Rademacher theorem \citep[Theorem 9.60]{rockafellar2009variational} and $T^{-1}$ is Lipschitz, we have
	\[
	m\left(y \in D^c \cap S_2\right) = m\big( (\lambda, \xi) \in T^{-1}( D^c \cap S_2)\big) = 0.
	\]
	Let $S_3\coloneqq\{b\in\mathbb{R}^d: \|b\|\leq 1, b^\top (x_{t,k}-x_t) = 0\}$.
By \citep[Corollary 4]{barthe2005probabilistic}, we have $b_{t,k+1} \sim \mathop{\textnormal{Unif}}(S_3)\overset{d}{=}X^\bot \xi$, where $\xi \sim \mathop{\textnormal{Unif}}(\mathbb{B}^{d-1})$. With $\lambda \sim \mathop{\textnormal{Unif}}([0,1])$ and countable union of zero measure set is negligible, we have
\[
\mathbb{P}\big(y_{t,k} \in D^c\cap S_2, \forall (t,k)\in[T]\times [K]\big) = 0,
\]
which completes the proof.
\end{proof}

\lemm*
\begin{proof}
	Let $\mathcal{F}_{t,k}=\sigma(y_{t,1},\cdots,y_{t,k})$ and $\widehat{\mathcal{F}}_{t,k}=\sigma(y_{t,1},\cdots,y_{t,k},b_{t,k+1})$. We denote $D_{t,k}$ as the event that $k$-loop does not break at $x_{t,k}$, i.e., $\|m_{t,k}\|>\epsilon$ and $f(x_{t,k})-f(x_t) > -\frac{\delta}{4}\|m_{t,k}\|$. It is clear that $D_{t,k}\in \mathcal{F}_{t,k}\subset \widehat{\mathcal{F}}_{t,k}$.
	Let $\gamma (\lambda) = (1-\lambda)x_t + \lambda \left(x_{t,k} + \frac{\delta\|m_{t,k}\|}{8L}\cdot b_{t,k+1}\right)$ for $\lambda \in [0,1]$. Note that $\gamma^\prime (\lambda) = x_{t,k}-x_t + \frac{\delta\|m_{t,k}\|}{8L}\cdot b_{t,k+1}$. Let $x_{t,k}^\prime = x_{t,k}+\frac{\delta\|m_{t,k}\|}{8L}\cdot b_{t,k+1}$. Since $y_{t,k+1}$ is uniformly sampled from the line segment $\left[x_t, x_{t,k}^\prime\right]$ and $f$ is differentiable at $y_{t,k+1}$ almost surely by \Cref{lem:zeroprob}, we know that
	\[
	\mathbb{E}\left[\left.\left\langle g_{t,k+1}, x_{t,k}^\prime-x_t  \right\rangle\right|\widehat{\mathcal{F}}_{t,k}\right]=\int_0^1 f^\prime (\gamma(t); x_{t,k}^\prime-x_t)\d t = f(x_{t,k}^\prime)-f(x_t).
	\]
	By $x_{t,k}^\prime-x_t=-\left(1-\frac{\|m_{t,k}\|}{8L}\right)\cdot\delta \frac{m_{t,k}}{\|m_{t,k}\|}+\frac{\delta\|m_{t,k}\|}{8L}\cdot b_{t,k+1}$, we have
	
	\[
	\begin{aligned}
	&\ \mathbb{E}\left[\left.\left\langle g_{t,k+1}, m_{t,k}  \right\rangle\right|\widehat{\mathcal{F}}_{t,k}\right]\\
	=&\ -\frac{\|m_{t,k}\|}{\left(1-\frac{\|m_{t,k}\|}{8L}\right)\cdot\delta}\cdot\mathbb{E}\left[\left.\left\langle g_{t,k+1}, x_{t,k}^\prime-x_t  \right\rangle\right|\widehat{\mathcal{F}}_{t,k}\right]+\frac{\|m_{t,k}\|}{\left(1-\frac{\|m_{t,k}\|}{8L}\right)\cdot\delta}\cdot \mathbb{E}\left[\left.\left\langle g_{t,k+1},  \frac{\delta\|m_{t,k}\|}{8L}\cdot b_{t,k+1} \right\rangle\right|\widehat{\mathcal{F}}_{t,k}\right]\\
	\leq&{}\ -\frac{\|m_{t,k}\|}{\left(1-\frac{\|m_{t,k}\|}{8L}\right)\cdot\delta}\cdot\Big(f(x_{t,k})-f(x_t) - | f(x_{t,k}^\prime) - f(x_{t,k})|\Big) + \frac{\|m_{t,k}\|^2}{8\left(1-\frac{\|m_{t,k}\|}{8L}\right)} \\
	\leq&{}\ -\frac{\|m_{t,k}\|}{\left(1-\frac{\|m_{t,k}\|}{8L}\right)\cdot\delta}\cdot\Big(f(x_{t,k})-f(x_t)  \Big) 
	+\frac{\|m_{t,k}\|^2}{4\left(1-\frac{\|m_{t,k}\|}{8L}\right)},
	\end{aligned}
	\]
	which directly implies
	\[
	\mathbb{E}\left[\left.\left\langle g_{t,k+1}, m_{t,k}  \right\rangle\right|\mathcal{F}_{t,k}\right] \leq 
	-\frac{\|m_{t,k}\|}{\left(1-\frac{\|m_{t,k}\|}{8L}\right)\cdot\delta}\cdot\Big(f(x_{t,k})-f(x_t)  \Big) 
	+\frac{\|m_{t,k}\|^2}{4\left(1-\frac{\|m_{t,k}\|}{8L}\right)}.
	\]
	By construction, $m_{t,k+1}=\beta m_{t,k} + (1-\beta) g_{t,k+1}$ under $D_{t,k}\cap \cdots \cap D_{t,1}$, and $m_{t,k+1}=0$ otherwise. Let $\mathcal{D}_{t,k}=D_{t,k}\cap \cdots \cap D_{t,1}$. Therefore,
	\begin{align*}
	&\ \mathbb{E}\left[\left. \|m_{t,k+1}\|^2 \right|\mathcal{F}_{t,k}\right] \\
	\leq &{}\ \left( \beta^2 \|m_{t,k}\|^2 + (1-\beta)^2L^2 +2\beta(1-\beta)\cdot \mathbb{E}\left[\left.\left\langle g_{t,k+1}, m_{t,k}  \right\rangle\right|\mathcal{F}_{t,k}\right] \right) \mathbb{1}_{\mathcal{D}_{t,k}} \\
	\leq &{}\ \left(\beta^2 \|m_{t,k}\|^2 + (1-\beta)^2L^2 +2\beta(1-\beta)\cdot \left(-\frac{\|m_{t,k}\|}{\left(1-\frac{\|m_{t,k}\|}{8L}\right)\cdot\delta}\cdot\Big(f(x_{t,k})-f(x_t)  \Big) 
	+\frac{\|m_{t,k}\|^2}{4\left(1-\frac{\|m_{t,k}\|}{8L}\right)}\right) \right)\mathbb{1}_{\mathcal{D}_{t,k}} \\
	\leq &{}\ \left(\beta^2 \|m_{t,k}\|^2 + (1-\beta)^2L^2 +2\beta(1-\beta)\cdot \frac{\|m_{t,k}\|^2}{2\left(1-\frac{\|m_{t,k}\|}{8L}\right)}\right)\mathbb{1}_{\mathcal{D}_{t,k}}\eqqcolon h(\beta) \mathbb{1}_{\mathcal{D}_{t,k}}.
	\end{align*}
By rearranging, it holds
\[
h(\beta) = \beta^2\underbrace{\left(L^2+\|m_{t,k}\|^2-\frac{\|m_{t,k}\|^2}{1-\frac{\|m_{t,k}\|}{8L}}\right)}_{T_1}  + \beta \underbrace{\left(\frac{\|m_{t,k}\|^2}{1-\frac{\|m_{t,k}\|}{8L}}-2L^2\right)}_{T_2} + L^2.
\]
Note that, by $0<\|m_{t,k}\|\leq L$ in $\mathcal{D}_{t,k}$ and conditioning on $\mathcal{D}_{t,k}$, it holds
\[
T_1 = \frac{L^2}{8\left(1-\frac{\|m_{t,k}\|}{8L}\right)}\left(8-\frac{\|m_{t,k}\|}{L}-\frac{\|m_{t,k}\|^3}{L^3}\right)>0, \quad\text{and}\quad T_2=\frac{L^2}{4\left(1-\frac{\|m_{t,k}\|}{8L}\right)}\left(-8+\frac{\|m_{t,k}\|}{L}+\frac{4\|m_{t,k}\|^2}{L^2}\right)<0.
\]
Thus, $h(\beta)$ achieves the minimum at $\beta_{t,k}=\frac{8L^3-L^2\|m_{t,k}\|-4L\|m_{t,k}\|^2}{8L^3-L^2\|m_{t,k}\|-\|m_{t,k}\|^3}$, which belongs to $\mathcal{F}_{t,k}$. Since $0<\|m_{t,k}\|\leq L$ in $\mathcal{D}_{t,k}$, we have
\[
h(\beta_{t,k}) =  \left(1-c_1 \cdot \frac{\|m_{t,k}\|^2}{L^2}\right)\|m_{t,k}\|^2 \overset{(\sharp)}{\leq} \left(1-\frac{\|m_{t,k}\|^2}{5L^2}\right)\|m_{t,k}\|^2,
\]
where
\[
c_1 = \frac{L^2}{\|m_{t,k}\|^2} - \frac{L^2h(\beta_{t,k})}{\|m_{t,k}\|^4}=\frac{16L^4 - 8L^3 \|m_{t,k}\| + L^2\|m_{t,k}\|^2}{64L^4-16L^3\|m_{t,k}\| + L^2\|m_{t,k}\|^2 - 8L\|m_{t,k}\|^3 + \|m_{t,k}\|^4}.
\]
Let $0<t\coloneqq\frac{\|m_{t,k}\|}{L}\leq 1$. For the validity of inequality $(\sharp)$, we observe
\[
\frac{1}{5}< \inf_{0<t\leq 1} \frac{16 - 8t + t^2}{64-16t + t^2 - 8t^3 + t^4}\leq c_1.
\]
To see it, note that, for $0<t\leq 1$, it holds
\[
5\cdot(16 - 8t + t^2) - (64-16t + t^2 - 8t^3 + t^4) = (t+2)\left(10\left(t-\frac{4}{5}\right)^2+\frac{8}{5}-t^3\right) > 0.
\]
Therefore,
	\[
	\mathbb{E}\left[\|m_{t,k+1}\|^2 \right]  
	=\mathbb{E}\Big[\mathbb{E}\left[\left. \|m_{t,k+1}\|^2 \right|\mathcal{F}_{t,k}\right] \Big] 
	\leq \mathbb{E}\Big[\left(1-\frac{\|m_{t,k}\|^2}{5L^2}\right)\|m_{t,k}\|^2 \Big]\leq \left(1-\frac{\mathbb{E}[\|m_{t,k}\|^2]}{5L^2}\right)\mathbb{E}[\|m_{t,k}\|^2].
	\]
	Then, by a similar argument in the proof of \citep[Lemma 13]{zhang2020complexity} we have $\mathbb{E}[\|m_{t,K}\|^2]\leq\frac{5L^2}{K+4}$.
	When $K\geq\frac{80L^2}{\epsilon^2}$, we have $\mathbb{E}[\|m_{t,K}\|^2]\leq \frac{\epsilon^2}{16}$.
\end{proof}

\thmconicdet*
\begin{proof}
Using Lemma~\ref{lem:zeroprob} and Lemma~\ref{lem:m}, the remaining parts directly follow from the proof of  \citep[Theorem 8]{zhang2020complexity}.	
\end{proof}
\section{Proofs of \Cref{sec:sto}}
\thmSINGD*
\begin{proof}
Let $\alpha_i \coloneqq \beta^{t-i}(1 - \beta)$ and denote $x'_{t+1} \coloneqq x_{t+1} + \zeta b_{t+1}$ and $\mathcal{G}_t \coloneqq \sigma (\sg_1, \ldots, \sg_t), \forall t$. Clearly, the random variables $m_t, x_t, x_{t+1}, \eta_t$ are $\mathcal{G}_t$-measurable. Note that
\[
m_t = \beta^{K} m_{t-K} + \sum_{i = t-K+1}^t {\alpha_i \sg_i}.
\]

Conceptually, if we choose $K$ to be sufficiently large, the term $\beta^{K} m_{t-K}$ is negligible. Then, if all the points $y_{t-K+1}, \ldots, y_t$ are inside $x_{t - K} + \delta \mathbb{B}$,  we have that $m_t$ approximately belongs to $\partial_\delta f(x_{t-K})$ in expectation. 

Note that for all $i = t-K+1,\ldots, t$,
\[
\begin{aligned}
	\norm{y_i - x_{t-K}} &\leq \norm{y_i - x_{i-1}} + \norm{x_{i-1} - x_{t-K}} \\
	&\mleq{a} \norm{x'_{i} - x_{i-1}} + \norm{x_{i-1} - x_{t-K}} \\
	&= \norm{\zeta b_i - \eta_{i-1} m_{i-1}} + \norm{\sum_{j = t-K}^{i-2} {\eta_j m_j}}\\
	&\mleq{b} \zeta + \sum_{j=t-K}^{i-1}{\eta_j\norm{m_j}} \\
	&\mleq{c} \frac{\omega}{p} + \frac{i-t+K}{p} \\
	&\leq \frac{K+\omega}{p},
\end{aligned}
\] 
where $\mar{a}$ holds since $y_i$ is sampled from the line segment $[x_{i-1}, x'_i]$, $\mar{b}$ uses $\norm{b_i} \leq 1$ and $\mar{c}$ follows from $\zeta \leq \frac{\omega}{p}$ and $\eta_t\norm{m_t} \leq \frac{1}{p}, \forall t$. We verify that the choices of $K, \omega$ and $p$ satisfy $\frac{K+\omega}{p} \leq \delta$:
\[
\frac{K+\omega}{p} = \frac{\frac{1}{1-\beta}\ln{\frac{16G}{\epsilon}}}{\frac{64G^2}{\delta\epsilon^2}\ln{\frac{16G}{\epsilon}}} = \delta.
\]
Then, conditioned on $\mathcal{G}_{t-K}$, since for all $i = t-K+1,\ldots, t$,
\[
\E{\sg_i \mid \mathcal{G}_{t-K}} = \E{\pf{y_i} \mid \mathcal{G}_{t-K}} \in \partial_\delta f(x_{t-K}),
\]
we have (note that $\sum_{i = t-K+1}^t {\alpha_i} = 1 - \beta^K$)
\[
\begin{aligned}
	&\ \ \ \ \frac{1}{1 - \beta^K}\sum_{i = t-K+1}^t {\alpha_i \E{\sg_i\mid \mathcal{G}_{t-K}}} \in \partial_\delta f(x_{t-K})\\
	&\Rightarrow \frac{1}{1 - \beta^K}\left(\E{m_t\mid \mathcal{G}_{t-K}} - \beta^K m_{t-K}\right) \in \partial_\delta f(x_{t-K}) \\
	&\begin{aligned}
		{}\Rightarrow \textup{dist}(0, \partial_\delta f(x_{t-K})) &\leq \frac{1}{1 - \beta^K}\left(\norm{\E{m_t\mid \mathcal{G}_{t-K}}} + \beta^K \norm{m_{t-K}}\right)\\ &\leq \frac{1}{1 - \beta^K}\left(\E{\norm{m_t}\mid \mathcal{G}_{t-K}} + \beta^K \norm{m_{t-K}}\right).
	\end{aligned}
\end{aligned}
\]

Take expectation.
\begin{gather}
	\E{\textup{dist}(0, \partial_\delta f(x_{t-K}))} \leq  \frac{1}{1 - \beta^K}\E{\norm{m_t}} + \frac{\beta^K G}{1 - \beta^K}, \nonumber\\
	\frac{1}{T}\sum_{t=1}^{T}{\E{\textup{dist}(0, \partial_\delta f(x_{t-K}))}} \leq  \frac{1}{(1 - \beta^K)T}\sum_{t=1}^T{\E{\norm{m_t}} }+ \frac{\beta^K G}{1 - \beta^K}.\nonumber
\end{gather}

We verify that the choices of $\beta$ and $K$ satisfy $\beta^KG \leq \frac{\epsilon}{16}$:
$\left(\beta^K  \leq \frac{\epsilon}{16G}\right) \Leftrightarrow \left(K \geq \frac{1}{\ln{\frac{1}{\beta}}}\ln{\frac{16G}{\epsilon}}\right)$. WLOG, we assume that $\epsilon \leq G$, and thus $\beta^K \leq \frac{1}{16}$. The above inequality can be further bounded as
\begin{equation}
	\label{ST3}
	\frac{1}{T}\sum_{t=1}^{T}{\E{\textup{dist}(0, \partial_\delta f(x_{t-K}))}} \leq  \frac{16}{15T}\sum_{t=1}^T{\E{\norm{m_t}} }+ \frac{\epsilon}{15}.
\end{equation}

The remaining proof is to show that Algorithm~\ref{alg:SINGD} ensures that $\frac{1}{T}\sum_{t=1}^T{\E{\norm{m_t}} } = O(\epsilon)$.

For ease of analysis, we denote $\mathcal{Y}_{t+1} \coloneqq \sigma (\sg_1, \ldots, \sg_t, b_{t+1}, y_{t+1})$ and $\hat{\mathcal{Y}}_{t+1} \coloneqq \sigma(\sg_1, \ldots, \sg_t, b_{t+1})$. Clearly, we have $\mathcal{G}_t \subset \hat{\mathcal{Y}}_{t+1} \subset \mathcal{Y}_{t+1} \subset \mathcal{G}_{t+1}$. Let $\varphi(\lambda) \coloneqq (1 - \lambda) x_t + \lambda x'_{t+1}$ for $\lambda \in [0,1]$. Since $y_{t+1}$ is uniformly sampled from the line segment $[x_t, x_{t+1}']$ and that $f$ is differentiable at $y_{t+1}$ almost surely, it holds that
\begin{equation}\label{ST1}
	\begin{aligned}
		\E{\innr{\sg_{t+1}, x'_{t+1} - x_t} \mid \mathcal{G}_t} &= \E{\E{\E{\innr{\sg_{t+1}, x'_{t+1} - x_t} \mid \mathcal{Y}_{t+1}} \mid \hat{\mathcal{Y}}_{t+1}} \mid \mathcal{G}_t} \\
		&= \E{\E{\innr{\pf{y_{t+1}}, x'_{t+1} - x_t}\mid \hat{\mathcal{Y}}_{t+1}} \mid \mathcal{G}_t} \\
		&= \E{\int_0^1{f'(\varphi(\lambda); x'_{t+1} - x_t) \textup{d}\lambda} \mid \mathcal{G}_t} \\
		&= \E{f(x'_{t+1}) - f(x_t) \mid \mathcal{G}_t}.
	\end{aligned}
\end{equation}

By $x'_{t+1} - x_t = - \eta_t m_t + \zeta b_{t+1}$, we have
\[
\begin{aligned}
	\E{\innr{\sg_{t+1}, x'_{t+1} - x_t} \mid \mathcal{G}_t} &= - \eta_t\E{\innr{\sg_{t+1},  m_t}\mid \mathcal{G}_t} + \zeta\E{\innr{\sg_{t+1}, b_{t+1}} \mid \mathcal{G}_t} \\
	&\leq - \eta_t\E{\innr{\sg_{t+1},  m_t}\mid \mathcal{G}_t} + \zeta G,
\end{aligned}
\]
where we used $\norm{b_{t+1}} \leq 1$. Thus, combining with \eqref{ST1}, we obtain
\begin{equation}\label{ST2}
	\begin{aligned}
		\E{\innr{\sg_{t+1},  m_t}\mid \mathcal{G}_t} &\leq \frac{1}{\eta_t} \E{f(x_t) - f(x_{t+1}) + f(x_{t+1}) - f(x'_{t+1}) \mid \mathcal{G}_t} + \frac{\zeta}{\eta_t} G \\
		&\leq \frac{1}{\eta_t} \big(f(x_t) - f(x_{t+1})\big) + \frac{\zeta}{\eta_t} (L + G).
	\end{aligned}
\end{equation}

Based on the construction $m_{t+1} = \beta m_t + (1 - \beta) \sg_{t+1}$, we can conclude that
\[
\begin{gathered}
	\norm{m_{t+1}}^2 = \beta^2\norm{m_t}^2 + 2\beta(1-\beta) \innr{\sg_{t+1}, m_t} + (1-\beta)^2 \norm{\sg_{t+1}}^2,
	\\
	\E{\eta_t\big(\norm{m_{t+1}}^2 - \beta^2\norm{m_t}^2 \big)} =  2\beta(1-\beta) \E{\eta_t\innr{\sg_{t+1}, m_t}} + (1-\beta)^2 \E{\eta_t\norm{\sg_{t+1}}^2}.
\end{gathered}
\]

From \eqref{ST2}, it holds that
\[
\begin{aligned}
	\E{\eta_t\big(\norm{m_{t+1}}^2 - \beta^2\norm{m_t}^2 \big)} \leq{}& 2\beta(1-\beta) \E{f(x_t) - f(x_{t+1})} + 2\beta(1-\beta) (L + G)\zeta \\ &+ (1-\beta)^2 \E{\eta_t\norm{\sg_{t+1}}^2},\\
	\frac{1}{T}\sum_{t=1}^T{\E{\eta_t\big(\norm{m_{t+1}}^2 - \beta^2\norm{m_t}^2 \big)}} \leq{}& \frac{2\beta(1-\beta)\Delta}{T} + 2\beta(1-\beta) (L + G) \zeta + \frac{(1-\beta)^2 G^2}{q},
\end{aligned}
\]
where we used $\eta_t \leq \frac{1}{q}$. 

Since $\eta_t = \frac{1}{p\norm{m_t} + q}$, using the same telescoping proof in \citep{zhang2020complexity}, as long as $\frac{pG}{q} \leq \frac{\beta}{2}$, the following holds
\[
\frac{1}{T}\sum_{t=1}^T{\E{\eta_t\big(\norm{m_{t+1}}^2 - \beta^2\norm{m_t}^2 \big)}} \geq \frac{\beta(1 - \beta)}{2T} \sum_{t = 1}^{T+1} {\E{\frac{\norm{m_t}^2}{p\norm{m_t} + q}}} - \frac{\beta G^2}{qT}.
\]

Thus, 
\[
\begin{gathered}
	\frac{\beta(1 - \beta)}{2T} \sum_{t = 1}^{T+1} {\E{\frac{\norm{m_t}^2}{p\norm{m_t} + q}}} \leq{} \frac{2\beta(1-\beta)\Delta}{T} + 2\beta(1-\beta) (L + G) \zeta + \frac{(1-\beta)^2 G^2}{q} + \frac{\beta G^2}{qT}, \\
	\frac{1}{T} \sum_{t = 1}^{T} {\E{\frac{q\norm{m_t}^2}{p\norm{m_t} + q}}} \leq{} \frac{4q\Delta}{T} + 4q (L + G) \zeta + \frac{2(1-\beta) G^2}{\beta} + \frac{2 G^2}{ T(1 - \beta)}.
\end{gathered}
\]

Comparing the above inequality with (14)\footnote{There is a typo in the telescoping proof of Theorem 14 in \citep{zhang2020complexity}: The term $\frac{\beta^2 G^2}{q}$ above Equation (14) should be $\frac{\beta G^2}{q}$. This typo does not affect the final convergence result.} in \citep{zhang2020complexity}, we notice that the only difference is the additional perturbation term $4q (L + G) \zeta$. Since we choose the identical $\beta, p, q$ and $T$ as in \citep{zhang2020complexity}, using the arguments (15) and (16) in \citep{zhang2020complexity} and denoting $m_{avg} \coloneqq \frac{1}{T}\sum_{t=1}^T{\E{\norm{m_t}} }$, we obtain
\[
\begin{aligned}
	\frac{4Gm_{avg}^2}{m_{avg} + 4G} &\leq \frac{\epsilon^2}{17} + 4q (L + G)\zeta \\
	&\mleq{\star} \frac{\epsilon^2}{15},
\end{aligned}
\]
where $\mar{\star}$ uses $\zeta \leq \frac{\epsilon^2}{510q(L+G)}$. The above is a quadratic equation in $m_{avg}$:
\[
	4Gm_{avg}^2 - \frac{\epsilon^2}{15} m_{avg} - \frac{4G\epsilon^2}{15} \leq 0.
\]
Solving for the positive root of this quadratic equation and using $\epsilon\leq G$, we obtain
\[
m_{avg} \leq \frac{\frac{\epsilon^2}{15} + \sqrt{\frac{\epsilon^4}{225} + \frac{64G^2\epsilon^2}{15}}}{8G} \leq \frac{4\epsilon}{15}  \leq \frac{5\epsilon}{16}.
\]

Finally, using \eqref{ST3}, we conclude that
\[
	\E{\textup{dist}(0, \partial_\delta f(x_{out}))} = \frac{1}{T}\sum_{t=1}^{T}{\E{\textup{dist}(0, \partial_\delta f(x_{t-K}))}} \leq \frac{2\epsilon}{5}.
\]
Thus, with probability at least $\frac{3}{5}$, we have $\textup{dist}(0, \partial_\delta f(x_{out})) \leq \epsilon$. 
\end{proof}

\section{Proofs of \Cref{sec:olcmain}}

\thmmain*
\begin{proof}
	As $x$ is Goldstein $(\epsilon,\eta)$-stationary, we have $\dist\big(0,\partial f(x+\eta\mathbb{B})\big)\leq \epsilon$, which implies that there exists 
	\[
	\|g\|\leq \epsilon, \qquad \textnormal{such that}\qquad g\in\partial f(x+\eta\mathbb{B})=\conv \left\{ \bigcup_{y\in\mathbb{B}_{\eta}(x)} \partial f(y) \right\}.
	\]
	By Carath\'eodory's theorem \citep[Theorem 2.29]{rockafellar2009variational}, we can write $g=\sum_{j=1}^{d+1} \alpha_j g_j$, where $\alpha_j \geq 0, \sum_{j=1}^{d+1}\alpha_j=1,g_j \in \partial f(y_j),y_j \in \mathbb{B}_{\eta}(x), \forall j \in [d+1]$.

	Let $y \in \mathbb{B}_\delta(x)$ be a pivot such that $\partial f$ is $\kappa$-outer Lipschitz continuous on $\mathbb{B}_\eta(x).$ As $f$ is Lipschitz and by \citep[Proposition 2.1.2]{clarke1990optimization}, $\partial f(y)$ is nonempty, convex, and compact.
	Let $u_j \coloneqq \arg \min_{z \in \partial f(y)} \|z-g_j\|, u\coloneqq\sum_{j=1}^{d+1}\alpha_j u_j \in \partial f(y)$. Then, we compute
	\[
	\|u\| = \left\| \sum_{j=1}^{d+1}\alpha_j u_j\right\| 
	\leq \|g\| + \sum_{j=1}^{d+1}\alpha_j \|u_j-g_j\|
	\leq \|g\| + \kappa\sum_{j=1}^{d+1}\alpha_j  \|y-y_j\|
	\leq \|g\| + \kappa\sum_{j=1}^{d+1}\alpha_j  \big(\|y-x\| + \|x - y_j\|\big)
	\leq \epsilon + \kappa (\delta+\eta),
	\]
	which completes the proof.
\end{proof}

\thmgen*
\begin{proof}
Let $\bigcup_{x \in S} \mathbb{B}^\circ_{\frac{1}{2}\bar{\delta}(x)} (x)$ be an open cover of $S$, where $\bar{\delta}(x)=\min\{\delta, \delta(x)\}$ and $\delta(x)$ is the inradius of neighborhood $V(x)$, on which $\partial f$ is $\kappa$-outer Lipschitz at $x$, satisfying $\mathbb{B}_{\delta(x)}(x) \subseteq V(x)$. As $S$ is compact, we find a finite subcover $\bigcup_{i \in [m]} \mathbb{B}^\circ_{\frac{1}{2}\delta_i} (x_i)$, where $\delta_i = \min\{\delta, \delta(x_i)\}$ and $x_i\in S$. Let $\eta \coloneqq \min_{i\in[m]} \frac{\delta_i}{2}$. Then, by Lebesgue's number theorem \citep[Chapter 3, Lemma 7.2]{munkres1974topology} on open cover $\bigcup_{x \in S} \mathbb{B}^\circ_{\bar{\delta}(x)} (x)$ of $S$, for any $x\in S$, there exists $i\in[m]$ such that $\mathbb{B}^\circ_{\eta}(x) \subseteq \mathbb{B}^\circ_{\delta_i}(x_i)$. Thus, $\mathbb{B}_{\eta}(x) \subseteq \mathbb{B}_{\delta_i}(x_i)\subseteq V(x_i)$. For any $z \in \mathbb{B}_\eta(x)\cap S$, by $\kappa$-outer Lipschitz continuity $\partial f$ on $V(x_i)$, we have
\[
\partial f(z) \subseteq \partial f(x_i) + \kappa \|x_i-z\| \mathbb{B},\ \forall z \in \mathbb{B}_\eta(x)\cap S,
\]
where $x_i \in \mathbb{B}_{\delta_i}(x)\subseteq \mathbb{B}_{\delta}(x)$. This completes the proof.
\end{proof}

\section{Proofs of \Cref{sec:rules}}
\propcaclsmooth*
\begin{proof}
	Let $F\coloneqq f+g$. By \citep[Exercise 8.8(c)]{rockafellar2009variational}, $\partial F = \partial f + \nabla g$.
	Let $y\in\mathbb{B}_\delta(x)$ be a pivot of $\partial f$. Then, for $\forall z \in \mathbb{B}_\eta(x)$, we compute 
\[
\partial F(z) = \partial f(z) + \nabla g(z) \subseteq \partial f(z) +\nabla g(y) + \beta \|z-y\|\mathbb{B}
\subseteq \partial F(y) + (\beta+\kappa) \|z-y\|\mathbb{B},
\]
which completes the proof.
\end{proof}

\propcalcsum*
\begin{proof}
	By \citep[Proposition 10.5]{rockafellar2009variational} and $f$ is Lipschitz, $\partial f = \bigoplus_{i=1}^m \partial f_i$.
	Let $y_i\in\mathbb{B}^{d_i}_{\delta_i}(x_i)$ be a pivot of $\partial f_i$. Also $y\coloneqq\bigoplus_{i=1}^m y_i$. Similarly, for any $z\in\mathbb{B}_\eta(x)$, it holds $z_i \in \mathbb{B}_{\eta}^{d_i}(x_i)\subseteq \mathbb{B}_{\eta_i}^{d_i}(x_i), \forall i \in [m]$. We compute 
	\[
	\partial f(z) = \bigoplus_{i=1}^m \partial f_i (z_i) \subseteq \bigoplus_{i=1}^m \Big(\partial f_i (y_i) + \kappa_i \|y_i - z_i\| \mathbb{B}^{d_i} \Big)
	\subseteq \partial f(y) + \kappa \|z-y\|\mathbb{B}^d,
	\]
	where $\|y-x\|^2 = \sum_{i=1}^m \|y_i - x_i\|^2 \leq \sum_{i=1}^m \delta_i^2 = \delta^2$. This completes the proof.
\end{proof}

\propcalclinear*
\begin{proof}
	Let $F(x) \coloneqq f(Ax)$. As $A$ is surjective, by \citep[Exercise 10.7]{rockafellar2009variational}, $\partial F(x) = A^\top \partial f(Ax)$. Let $q \in \mathbb{B}^n_\delta(Ax)$ be a pivot of $\partial f$. Let $y\coloneqq A^\dag q + (I-A^\dag A)x$. Then $Ay=q$ and $\|y-x\|\leq \|A^\dag\| \|q - Ax\| \leq \delta\|A^\dag \|  $. Meanwhile, for any $z \in \mathbb{B}^d_{\frac{\eta}{\|A\|}}(x)$, it holds $\|Az - Ax\| \leq \|A\| \|z-x\| \leq  \eta$. We compute
	\[
	\partial F(z) = A^\top \partial f(Az) \subseteq A^\top \partial f(Ay) + \kappa \|Ay-Az\|A^\top \mathbb{B}^n
	\subseteq  \partial F(y) + \kappa \|A\|^2 \|y-z\|\mathbb{B}^d,
	\]
		which completes the proof.
\end{proof}

\propcalcrescal*
\begin{proof}
		Let $F=g\circ f$. By \citep[Theorem 2.3.9(ii)]{clarke1990optimization}, $\partial F(x) = \nabla g(f(x)) \cdot  \partial f(x)$.
	Let $y\in\mathbb{B}_\delta(x)$ be a pivot of $\partial f$. Then, for $\forall z \in \mathbb{B}_\eta(x)$, we compute 
	\[
	\begin{aligned}
	\partial F(z) = \nabla g(f(z)) \cdot  \partial f(z) &\subseteq \nabla g(f(z)) \cdot  \partial f(y) + \nabla g(f(z)) \cdot\kappa \|z-y\|\mathbb{B} \\
	&\subseteq \nabla g(f(y)) \cdot  \partial f(y) + (\beta L_1 + \kappa L_2) \|z-y\|\mathbb{B} \\
	&= \partial F(y) + (\beta L_1 + \kappa L_2) \|z-y\|\mathbb{B},
	\end{aligned}
	\]
	which completes the proof.
\end{proof}

\propcalcsumshared*
\begin{proof}
	Let $y = P(x)\in\mathbb{B}_{\delta_i}(x)$ be a pivot of $G_i(x)$, which by pivot sharing assumption should hold for all $i\in[m]$. Thus $\|y-x\|\leq \min_{i\in[m]}\delta_i =  \delta$. For all $z\in\mathbb{B}_\eta(x)\subseteq \mathbb{B}_{\eta_i}(x)$, we compute 
	\[
	G(z) = \sum_{i=1}^m G_i(z) \subseteq \sum_{i=1}^m\Big(G_1(y) + \kappa_i \|z- y\|\mathbb{B}\Big)	\subseteq G(y) + \left(\sum_{i=1}^m \kappa_i\right)\|z - y\| \mathbb{B},
	\]
	as expected.
\end{proof}

\propcalcpartial*
\begin{proof}
	Let $(y_0,y_i)\in\mathbb{B}^{d_0+d_i}_{\delta_i}\big((x_0, x_i)\big)$ be a pivot of $G_i$. Also $y\coloneqq\bigoplus_{i=0}^m y_i$. Similarly, for any $z\in\mathbb{B}_\eta(x)$, it holds $(z_0, z_i) \in \mathbb{B}_{\eta}^{d_0+d_i}\big((x_0, x_i)\big)\subseteq \mathbb{B}_{\eta_i}^{d_0+d_i}\big((x_0, x_i)\big), \forall i \in [m]$. We compute 
	\[
	G(z) = \sum_{i=1}^m G_i \big((z_0, z_i)\big) \subseteq \sum_{i=1}^m \Big(G_i \big((y_0, y_i)\big) + \kappa_i \big\|(z_0, z_i) - (y_0, y_i)\big\| \mathbb{B}^{d_0+d_i} \Big)
	\subseteq G(y) + \kappa \|z-y\|\mathbb{B}^d,
	\]
	where $\|y-x\|^2 = \sum_{i=0}^m \|y_i - x_i\|^2 \leq \sum_{i=1}^m\big( \|y_0 - x_0\|^2 + \|y_i - x_i\|^2\big) \leq \sum_{i=1}^m \delta_i^2 = \delta^2,$ and $d = \sum_{i=0}^m d_i$. This completes the proof.
\end{proof}

\section{Proofs of \Cref{sec:application}}

\clmsubdrho*
\begin{proof}
	Define
	\begin{align*}
	C_1 &\coloneqq \{ (u_1, u_2) : u_2 \geq 0\}, \\
	C_2 &\coloneqq \{ (u_1, u_2) : u_2 \leq 0\}.
	\end{align*}
It is clear that $C_1 \cup C_2 = \mathbb{R}^2$, and we have
\[
\varrho(u_1, u_2) = 
\left\{ \begin{array}{ccl}
	u_1 \cdot u_2 & \mbox{for} 
	& (u_1, u_2) \in C_1, \\ 
	0 & \mbox{for} & (u_1, u_2) \in C_2.
\end{array}\right.
\]
Note that $C_1 \cap C_2$ form a set $S$ of measure $0$, and if $(u_1, u_2) \notin S$, then $\varrho$ is differentiable. The claim follows from taking convex hull with \citep[Theorem 9.61]{rockafellar2009variational}.
\end{proof}

In the following proof, we will use the following notion named partial Clarke subdifferential. See also \citep[Page 48]{clarke1990optimization}, \citep[Corollary 10.11]{rockafellar2009variational}.
\begin{definition}
	Let a local Lipschitz function $f:\mathbb{R}^n\times \mathbb{R}^m \rightarrow \mathbb{R}$ and $g_y: x\rightarrow f(x,y)$. Then the partial Clarke subdifferential with respect to the first argument is defined as $\partial_1 f(x,y) \coloneqq \partial g_y(x)$. $\partial_2 f(x,y)$ is defined similarly.
\end{definition}

\begin{claim}\label{clm:rhosubsep}
	$\partial \varrho (u_1, u_2) = \partial_1 \varrho(u_1, u_2) \times \partial_2 \varrho(u_1, u_2)$ and $|\pi_1 \circ \partial \varrho(u_1, u_2)| = 1$.
\end{claim}
\begin{proof}
	Note that $\partial_1 \varrho(u_1, u_2) = \max \{u_2, 0\}$ and $\partial_2 \varrho(u_1, u_2) = u_1 \cdot \partial(\max \{\cdot, 0\})(u_2)$. The proof completes by using \Cref{clm:subd-rho} and literally checking definitions.
\end{proof}

\clmchain*
\begin{proof}
To begin, we observe the following general fact.
For any set $A\subseteq \mathbb{R}^n\times \mathbb{R}^m$, if $|\pi_1 A|=1$, then $A = \pi_1 A \times \pi_2 A$. To see it, for one direction, if $(a_1, a_2) \in A$, then $a_1 \in \pi_1 A, a_2 \in \pi_2 A$. Thus, $A\subseteq \pi_1 A \times \pi_2 A$. For the other direction, let $a_1 \in \pi_1 A, a_2 \in \pi_2 A$. As $\{a_1\} = \pi_1 A$, then by the definition of $\pi_2 A$, it holds $(a_1, a_2) \in A$. Thus $ \pi_1 A \times \pi_2 A \subseteq A$.

To avoid uninformative sophisticated notation, we will use ``$\overset{P}{=}$'' for equivalence up to coordinate permutation. Formally, if $A\overset{P}{=}B$, then there exists a permutation matrix $P$ such that $B=\{P x: x \in A\}$. We compute
	\begin{align*}
	\partial f(a,U) &\subseteq \sum_{i=1}^n \partial h_i (a, u_i) \tag{\citet[Proposition 2.3.3]{clarke1990optimization}} \\
	&= \sum_{i=1}^n \nabla \ell_i \left( \sum_{j=1}^m \varrho(a_j, u_{ij}) \right)\cdot \bigoplus_{j=1}^m \partial \varrho(a_j, u_{ij}) \tag{\citet[Theorem 2.3.9(ii)]{clarke1990optimization}} \\
	&= \sum_{i=1}^n \nabla \ell_i \left( \sum_{j=1}^m \varrho(a_j, u_{ij}) \right)\cdot \bigoplus_{j=1}^m \Big( \partial_1 \varrho(a_j, u_{ij}) \times \partial_2 \varrho(a_j, u_{ij}) \Big) \tag{\Cref{clm:rhosubsep}}\\
	&\overset{P}{=}\underbrace{\left( \sum_{i=1}^n \nabla \ell_i \left( \sum_{j=1}^m \varrho(a_j, u_{ij}) \right)\cdot \bigoplus_{j=1}^m \partial_1 \varrho(a_j, u_{ij}) \right)}_{S_1}\times \underbrace{ \left( \bigoplus_{i=1}^n \bigoplus_{j=1}^m \nabla \ell_i \left( \sum_{j=1}^m \varrho(a_j, u_{ij}) \right)\cdot\partial_2 \varrho(a_j, u_{ij}) \right)}_{S_2}.
	\end{align*} 

	Note that $|S_1|=1 $ as $|\partial_1 \varrho (a_j, u_{ij})|=1, \forall (i,j) \in [n]\times [m]$ by \Cref{clm:rhosubsep}. Thus $|\pi_1 \circ \partial f(a,U)|=1$ and $\partial f(a,U) \overset{P}{=} \pi_1 \circ \partial f(a,U) \times \pi_2 \circ \partial f(a,U)$.
	With \citep[Proposition 2.3.16]{clarke1990optimization}, we compute
	\[
	\begin{aligned}
	\partial f(a,U) &\overset{P}{=} \Big(\pi_1 \circ \partial f(a,U)\Big) \times \Big(\pi_2 \circ \partial f(a,U)\Big) \\ 
	&\supseteq \partial_1 f(a,U) \times \partial_2 f(a,U) \\
	&\overset{\natural}{=} S_1 \times S_2.
	\end{aligned} 
	\]
	To see $(\natural)$, note that $f(\cdot, U)$ is differentiable. Thus it is straightforward to check $S_1 = \partial_1 f(a,U)$. For $S_2$, note that $f(a, \cdot)$ is fully separable (as $a$ is fixed). Then, with \citep[Proposition 10.5]{rockafellar2009variational} and \citep[Theorem 2.3.9(ii)]{clarke1990optimization}, the verification of $S_2 = \partial_2 f(a,U)$ is routine.
	
	This completes the proof.
\end{proof}
\end{document}